\documentclass[12pt]{article}
\usepackage{latexsym}
\usepackage{amssymb}
\usepackage{amsmath}
\usepackage{amsthm}
\usepackage{cases}
\def\UNO{1\mkern-7mu1}

\numberwithin{equation}{section}

\font\gorditas = msbm8

\def\zetita{\hbox{\gorditas Z}}
\def\Z{\hbox{\gorditas Z}}

\theoremstyle{plain}
\newtheorem{theorem}{Theorem}[section]

\newtheorem{proposition}[theorem]{Proposition}
\theoremstyle{definition}

\newtheorem{remark}[theorem]{Remark}

\begin{document}

\noindent

\begin{center}
{\bf\large  Oscillatory Fractional Brownian Motion \\[.4cm]
and Hierarchical Random Walks$^*$}\\[1cm]
\begin{center}
\begin{tabular}{lll}
Tomasz Bojdecki$^1$&Luis G. Gorostiza$^\dagger$&Anna Talarczyk$^1$\\
tobojd@mimuw.pl.edu&lgorosti@math.cinvestav.mx&annatal@mimuw.pl.edu
\end{tabular}
\end{center}

\end{center}
\vglue1cm
\centerline{\bf Abstract}
\vglue.5cm
We introduce oscillatory analogues of fractional Brownian 
motion, sub-fractional Brownian motion and other related long range dependent 
Gaussian processes, we discuss their properties, and we show how they arise 
from particle systems with or without branching and with different types 
of initial conditions, where the individual particle motion is the 
so-called $c$-random walk on a hierarchical group. 
The oscillations
are caused by the discrete and ultrametric structure of the hierarchical 
group, and they become slower as time tends to infinity and faster as time approaches zero. We also give other results 
to provide an overall picture of the behavior of this kind of systems, 
emphasizing the new phenomena that are caused by the ultrametric structure as compared with results for analogous models on Euclidean space.
\vglue.5cm
\noindent
\vglue.5cm
\noindent

\footnote{\kern-.6cm $^*$ Supported in part by CONACyT Grant 98998 (Mexico) and MNiSzW Grant N N201 397537 (Poland).\\
$^1$ Institute of Mathematics, University of Warsaw, ul. Banacha, 02-097, Warszawa, Poland.\\
$^\dagger$ Centro de Investigaci\'on y de Estudios Avanzados, A.P. 14-740, M\'exico 07000 D.F., Mexico.}

\section{Introduction}

The main objectives of  this paper are to introduce oscillatory analogues of fractional Brownian motion, sub-fractional Brownian motion and other related long range dependent Gaussian processes, and study their properties, and to show 
how those processes arise from systems of hierarchical random walks. In this Introduction we give some relevant background and put the paper in context.

It is well known that {\it fractional Brownian motion  (fBm)} with Hurst parameter $H\in (0,1/2)\cup (1/2,1)$, whose covariance function is
\begin{equation}
\label{eq:1.1}
\frac{1}{2}(s^{2H}+t^{2H}-|s-t|^{2H}),\quad s,t\geq 0,
\end{equation}
has an intrinsic mathematical interest and a wide variety of areas of application that are often mentioned in papers on the subject (see \cite{ST} for basic background and \cite{DOT} for some applications). The main properties of fBm are that it is a continuous centered Gaussian process which is self-similar, with stationary increments and long range dependence, and it does not have the Markov property and is not a semimartingale (the case $H=1/2$ corresponds to ordinary Brownian motion). A closely related process is {\it sub-fractional Brownian motion} ({\it sfBm}) with parameter $H\in (0,1/2)\cup (1/2,1)$, whose covariance function is
\begin{equation}
\label{eq:1.2}
s^{2H}+t^{2H}-\frac{1}{2}[(s+t)^{2H}+|s-t|^{2H}],\quad s,t\geq 0.
\end{equation}
It shares the properties of fBm except  increment stationarity, and it has a faster decay of long range dependence (the case $H=1/2$ also corresponds to ordinary Brownian motion). This process was analyzed in \cite{BGT0} and in \cite{BGT1},
where it was shown that it arises from occupation time fluctuations of a branching system of
 particles moving according to a symmetric $\alpha$-stable
 L\'evy process in $\mathbb R^d$, and $H$ is restricted to $(1/2,1)$.
A similar result in \cite{BGT1} 
for a system without branching leads to fBm with $H$ restricted to  $(1/2, 3/4]$.
SfBm also appeared independently at the same time in \cite{DZ} in a different context as the even part of fBm (defined on the whole real line). In \cite{BGT9} it was found, again for a system of $\alpha$-stable processes, that sfBm also arises in a  natural and intrinsic way (in fact, more natural than fBm) not necessarily  related with branching. A different type of particle picture interpretation for fBm and sfBm with values of $H$ in the whole interval $(0,1)$, also based on a system of $\alpha$-stable processes, was given in \cite{BT}, where the role of time is played by a spatial parameter. In a series of papers
\cite{BGT0,BGT1,BGT2,BGT3,BGT4,BGT5,BGT6,BGT7,BGT8,BGT9}, dealing with various models within a large class, we have studied occupation time fluctuations of $\alpha$-stable particle systems that lead in particular to fBm, sfBm and other related processes. Some of those papers will be 
referred to below for comparisons or for some technical points. Other authors have recently investigated some analogous  models where those processes and  similar ones also arise \cite{BB,BZ,LX,M1,M2,M3,Sw}, further properties of sfBm and related stochastic calculus have been discussed \cite{BB,EN,ET,LY,LYPW,
Me,N,RT,S,SC,SY1,SY2,SZ,T1,T2,T3,T4,T5,YS,YSH}, 
a strong approximation has been devised \cite{GGL}, and an application 
has been proposed \cite{LLY}.
Another  process that possesses the main properties of fBm except increment stationarity, and has also been studied by several authors recently, is bi-fractional Brownian motion, which was introduced in \cite{HV}.
According to some authors, processes that do not have stationary increments may be more useful for modelling in some applications (e.g., \cite{BMG} 
in financial markets). 

The oscillatory analogues of fBm and sfBm are (centered Gaussian) {\it oscillatory fractional Brownian motion (ofBm)} $\xi^H$, with parameter $H\in (1/2,1)$ and covariance function
\begin{equation}
\label{eq:1.3}
E\xi^H_s\xi^H_t=
\frac{1}{2}(s^{2H}h_s+t^{2H}h_t-|s-t|^{2H}h_{|s-t|}),\quad s,t\geq 0, 
\end{equation}
and {\it oscillatory sub-fractional Brownian motion (osfBm)} $\zeta^H$, with parameter $H\in (1/2,1)$ and covariance function
\begin{equation}
\label{eq:1.4}
E\zeta^H_s\zeta^H_t=
s^{2H}h_s+t^{2H}h_t-\frac{1}{2}[(s+t)^{2H}h_{s+t}+|s-t|^{2H}h_{|s-t|}],\quad s,t\geq 0, 
\end{equation}
where $(h_t)_{t>0}$ is a periodic function in logarithmic scale, i.e., $h_t=h_{at}, t>0$, $a$ being a constant, $0<a<1$. This function oscillates between two positive values, slower as $t\to\infty$, faster as $t\to 0$. A detailed description of these processes and other related ones is given in the next section. If the function $h$ were constant, ofBm and osfBm would reduce to fBm and sfBm, respectively. The function $h$ that arises from the models cannot be constant, hence ofBm and osfBm may not be considered as extensions of fBm and sfBm, but they can be regarded as oscillatory analogues of them.

Another process is  {\it oscillatory negative sub-fractional Brownian motion} ({\it onsfBm}), which has a covariance of the same form as that of osfBm, but with $H\in (1,3/2)$ and a slightly different oscillatory and negative function instead of $h$. This is an oscillatory analogue of the {\it negative sub-fractional Brownian motion (nsfBm)} that was found in \cite{BGT7} with a high density limit and $H\in (1,5/4)$. In \cite{BT} there is also a different type of
  particle picture interpretation for nsfBm with the full range $H\in (1,2)$.

It could be that the oscillatory analogues of fBm, sfBm and nsfBm turn out to be useful for some applications due to the additional ``volatility'' introduced by the oscillations of their covariances and their paths, specially for small times, while 
keeping most of the other properties of the respective non-oscillatory processes.

The oscillatory analogues of fBm, sfBm and nsfBm come about when instead of systems of $\alpha$-stable processes we consider systems of (a class of) hierarchical random walks.

The state space for hierarchical random walks is the {\it hierarchical group of order $M$} (integer $>1$), also called {\it hierarchical lattice} in 
physics, which is defined by
\begin{equation}
\label{eq:1.5}
\Omega_M=\{x=(x_1,x_2,...): x_i\in \{0,1,\ldots, M-1\},\quad 
\Sigma_ix_i<\infty\},
\end{equation}
with addition componentwise mod $M$. It is a countable Abelian group which is the direct sum of a countable number of copies of the cyclic group of order $M$. A translation invariant metric $|\cdot|$ on $\Omega_M$, called {\it hierarchical distance}, is defined by
\begin{equation}
\label{eq:1.6}
|x-y|=\left\{\begin{array}{lll}
0&{\rm if}&x=y,\\
\max\{i:x_i\neq y_i\}&{\rm if}&x\neq y,
\end{array}\right.
\end{equation}
and  satisfies the strong (non-Archimedean) triangle inequality
$$|x-y|\leq \max \{|x-z|, |z-y|\},\quad x,y,z\in \Omega_M.$$
Therefore $(\Omega_M,|\cdot|)$ is an {\it ultrametric space}. Background on ultrametric spaces  can be found e.g. in \cite{Sc}. The main difference with spaces of Euclidean type is that in an ultrametric space two balls are either disjoint or one is contained in the other. As a consequence, contrary to the Euclidean case, the only way to ``go out'' of a ball of radius $R$ is to make a single jump of size larger than $R$, and it is not possible by a sequence of small steps. $(\Omega_M,|\cdot |)$ can be represented as the set of leaves at the top of an infinite regular tree where each inner node at each level $j\geq 1$ from the top has one predecessor at level $j+1$ and $M$ successors at level $j-1$. The distance between two individuals (leaves) at the top is the depth in the tree to their most recent common node. Hierarchical structures of this type arise in the physical, biological, social and computer sciences due to the multiscale organization of many natural objects. For 
example, Sawyer and Felsenstein \cite{SF} have used a model like this where the distance represents the degree of genetic relatedness between individuals in a population that are grouped at several different levels of organization. Hierarchical structures have also been used in statistical physics, fundamentally in the work of Dyson \cite{D} on ferromagnetic models (see also 
\cite{CE,Si}), as they provide ``caricatures'' of Euclidean lattices which make it possible to carry out renormalization group analysis in a rigorous way. Some  papers that use hierarchical structures in connection with population models,
 branching systems, 
interacting diffusions, self-avoiding random walks, contact 
processes, percolation, and search algorithms are 
\cite{AS, BI, CDG,DGo1,DGo2,DGW1,DGW3,DG1,DG2,DGZ,FG,G,Kl,K}. Hierarchical random walks play a key role in some of the models in those papers, as do simple symmetric random walks in models on Euclidean lattices. 

A {\it hierarchical random walk} on $\Omega_M$ is defined by a probability distribution $(r_j)_{j\geq 1}$ on $\{1,2,\ldots\}$. The walk first chooses a distance $j$ with probability $r_j$, and then jumps to a point  
uniformly among all those that are at distance $j$ from its previous position. 
We assume that there exist arbitrarily large $j$ such that $r_j>0$, for otherwise the walk would never leave a bounded set, due to the abovementioned fact on 
jumping out of a ball. A particular hierarchical random walk that will be the most relevant in this paper is the {\it $c$-random walk}, where $c$ is a constant, $0<c<M$, and
\begin{equation}
\label{eq:1.7}
r_j=\left(1-\frac{c}{M}\right)\left(\frac{c}{M}\right)^{j-1},\quad j=1,2,\ldots
\end{equation}
A continuous time $c$-random walk with jump rate 1 has transition probability for going from $0$ to $y$ in time $t>0$ given  by
\begin{equation}
\label{eq:1.7a}
p_t(0,y)=e^{-t}\sum^\infty_{n=0}\frac{t^n}{n!}p^{(n)}(0,y),
\end{equation}
where $p^{(n)}(0,y)$ is the transition probability from $0$ to $y$ in $n$ steps for the discrete time walk. In what follows, $c$-random walk will always refer to the continuous time process.

The $c$-random walk ``mimics'' (or is a ``caricature'' of) a symmetric $\alpha$-stable process on $\mathbb R^d$ in a sense explained below. Some properties of $\alpha$-stable processes and $c$ random walks (and other L\'evy processes on additive Abelian groups) can be described in a unified way by means of the {\it degree} $\gamma$, which is the supremum over all the $\sigma >-1$ such that the operator power $G^{\sigma+1}$ of the potential operator $G$ of the underlying motion process is finite. For $\gamma>0$, the degree also gives information on how fast the process escapes to infinity. The degree is useful for comparison of results on particle systems involving $\alpha$-stable process and corresponding ones with $c$-random walks. The development of hierarchical random walks has been outlined in \cite{DGW2}, and their degree and several of their properties, in particular recurrence/transience, have been studied in \cite{DGW4} (see also a summary in the Appendix of \cite{BGT10}).

For symmetric $\alpha$-stable process on $\mathbb R^d, \gamma$ is given by
\begin{equation}
\label{eq:1.8}
\gamma=\frac{d}{\alpha}-1,
\end{equation}
and for $c$-random walk on $\Omega_M$ (either in discrete or continuous time), $\gamma$ is given by
\begin{equation}
\label{eq:1.9}
\gamma=\frac{\log c}{\log ({M}/{c})}.
\end{equation}
For $\alpha$-stable process $\gamma$ is restricted to the interval 
$[-{1}/{2},\infty)$, and for $c$-random walk $\gamma$ can take any value in $(-1,\infty)$. In this sense $c$-random walks are a richer class of processes. This is the reason why in the occupation time fluctuation limit referred to above for  $\alpha$-stable process without branching  the parameter $H$ (which is expressed in terms of $\gamma$ in the same way for both processes) in (\ref{eq:1.1}) is restricted to the interval $(1/2, 3/4]$, whereas for the $c$-random walk it can take any value in $(1/2,1)$, as we shall see in Theorem 2.2 and Remark 2.3(a). A similar thing happens with the nsfBm (for $\alpha$-stable process), which has $H\in (1,5/4)$ \cite{BGT7}, and the onsfBm (for $c$-random walk), which has $H\in (1,3/2)$ (Theorem 2.6). 

The sense in which $c$-random walks are considered to ``mimic'' $\alpha$-stable processes is that they have the same recurrence/transience behaviors for equal values of their degrees, and that they have the same spatial asymptotic decay of powers of their potential operators in terms of the ``Euclidean radial distance'' 
(see \cite{DGW4}, Remark 3.2.2(b)). But there are also differences, for example, the distance of the $c$-random walk from a given point behaves differently 
than the $\alpha$-stable Bessel process (see \cite{DGW4}, Remark 3.5.7). Branching  $c$-random walks can be considered similarly as branching $\alpha$-stable processes, 
with a given branching rate.

In \cite{DGW1}, occupation time fluctuation limits for branching particle systems were studied (with one and two branching levels, and also without branching), where the individual particle motion was either symmetric $\alpha$-stable process on $\mathbb R^d$ or $c$-random walk on $\Omega_M$. In both cases the results referred to the occupation time fluctuation at a final time $t$ as $t\to\infty$. A much more comprehensive asymptotics of the occupation time of the systems consists in accelerating the time $t$ by a scaling parameter $T$ and studying the limit of the fluctuation process for $t\in [0,\infty)$ as $T\to\infty$, as this reveals in detail the time behavior properties of the limit process. That was the objective of our papers in the case of $\alpha$-stable process on $\mathbb R^d$ (we considered the models without branching and with only one branching level, but also with other types of initial conditions and other differences with respect to 
\cite{DGW1}). 
The typical situation is that there are two different limiting regimes, corresponding to ``low'' and ``high'' dimensions $d$, and a ``critical'' dimension that represents the phase transition between the two regimes. Those dimensions depend on whether or not branching is present. For ``low'' dimensions the limit process has long range time dependence (fBm in the case without branching, and sfBm 
 with branching) and a uniform spatial structure (Lebesgue measure), for ``high'' dimensions it has independent increments in time and a non-trivial spatial structure (e.g., a non-trivial ${\cal S}'(\mathbb R^d)$-valued Wiener process), and at the ``critical'' dimension the fluctuations have a larger order of magnitude (with a logarithmic factor), the spatial structure is uniform and the time behavior is Brownian motion. Phase transition phenomena  are well known in several systems of interest in statistical physics.

Concerning  the spatial structure of the limits, the three behaviors recalled above were already found in \cite{DGW1} for both $\alpha$-stable process and $c$-random walk, and they were given in terms of $\gamma$ in a unified way. Namely, the regimes for the non-branching model in the ``low'', ``critical'' and ``high'' cases are $-1<\gamma<0, \gamma=0, \gamma>0$ (i.e., $1=d<\alpha, d=\alpha, d>\alpha; c<1, c=1, c>1$, respectively), and for the branching model (1 level branching), they are $0<\gamma<1, \gamma=1, \gamma>1$ (i.e, $\alpha< d<2\alpha, d=2\alpha, d> 2\alpha; 1<c<\sqrt{M}, c=\sqrt{M}, c>
\sqrt{M}$, respectively).

The natural question that remained was the possibility of extending results of 
\cite{DGW1} on $c$-random walks by adding the time structure, as described above,  analogously as it was done in our previous papers for $\alpha$-stable processes. This is the problem solved in the present paper. Clearly, the same  regimes occur with the phase transition between long range dependence in the ``low'' case and independent increments in the ``high'' case. The main novelty now is the appearance of 
the oscillatory analogues of fBm and sfBm
in the cases of long range dependence, namely, ofBm for $-1<\gamma<0$ (without branching), and osfBm for $0<\gamma <1$ (with branching), in both cases with Poisson initial condition, the latter also arising from a system without branching but with a different type of initial condition, and also  onsfBm (with branching and high density, $-1<\gamma<0$) 
(Theorems 2.2, 2.4 and 2.6). Generally, for the models with $\alpha$-stable processes or $c$-random walks, in the ``low'' cases the most interesting are the time properties of the limit processes, and in the ``high'' cases the spatial properties are the most relevant. 
Nevertheless, in all cases there are also significant differences between the two models. In the ``low'' cases the discrete and ultrametric structure of the state space $\Omega_M$ causes the time oscillations of the limit processes, and the limits exist only for special sequences $T_n\to\infty$. In the ``high'' cases, at a fixed time the spatial structure of the limit exhibits fractional behavior for $\alpha$-stable processes (see \cite{BT}), while for $c$-random walks it is Markovian of diffusion type (see Subsection 2.3.11).

In Section 2 we begin by describing the  systems of hierarchical random walks (interpreted as particle motions)
with different types of initial conditions that have effects on the limit processes, and we give the definitions of ofBm, osfBm, onsfBm, and other related processes that will appear. We then formulate convergence theorems for non-branching and branching systems, and discuss the results.
After that, we study properties of the obtained limit processes, including certain forms of self-similarity, long range dependence, non-Markov and non-semimartingale properties for ofBm and osfBm, and these two processes are represented in terms of a sequence of independent Brownian motions (Proposition 2.14). On the other hand, onsfBm is a semimartingale (by Proposition 2.16), analogously as nsfBm \cite{BGT5}.
We also give a result on quadratic variation (Proposition \ref{p:2.17a}) that  exhibits an oscillatory 
(in a sense)
character of paths of ofBm and osfBm.
Finally we point out the Markov property of the spatial structure in the ``high'' cases.

In Section 3 we give the proofs.
Since our main interest is identifying the limit processes, 
we prove functional convergence (i.e., tightness in addition to convergence of finite dimensional distributions) only for some of the main limits, focussing on the cases where the oscillatory limits occur. A new way of proving convergence that was introduced in \cite{BGT9}  (different from the method in our previous papers) is again useful here and allows to give reasonably short proofs of the limit theorems; some of the  ideas did not appear in the published version 
\cite{BGT9}, but they can be found in the complete version in arXiv. The methods of proof have some analogies to those used in the Euclidean case with stable processes, but the ultrametric structure of $\Omega_M$ involves different types of calculations.

In Section 4 several additional questions  are mentioned in a series of comments. They include results that complete the overall picture,  some conjectures
and some additional things that could be done with 
other models, but we do not go into  them because  they would not be significantly new and they would extend the length of the paper considerably. 

The following notation is used in the paper.

\parindent=0pt
${\cal B}_b$: space of real functions on $\Omega_M$ with bounded support;

${\cal B}'_b$: dual of ${\cal B}_b$;

$\langle\cdot,\cdot\rangle$: duality on ${\cal B}_b'\times {\cal B}_b$;

${\cal T}_t$: semigroup associated with $c$-random walk, i.e., 
\begin{equation}
\label{eq:1.10}
{\cal T}_t\varphi(x)=\sum_{y\in\Omega_M}p_t(x,y)\varphi(y), \,\,\varphi\in{\cal B}_b,
\end{equation}
where $p_t(x,y)=p_t(0,y-x)$ is given by (\ref{eq:1.7a});

$G$: potential operator (for $\gamma>0$), i.e., $G\varphi(x)= \int_0^\infty{\cal T}_t\varphi(x)dt$;

$\Rightarrow_{\kern-.35cm _f}\,\,$: weak convergence of finite dimensional distributions;

$\Rightarrow_{\kern-.35cm _c}\,\,$: weak convergence of distributions in the space of continuous functions  $C([0,\tau])$ for each $\tau>0$.

Generic constants are written $C$, $C_i$, with possible dependencies in parenthesis or in subscripts.
\parindent=20pt

\section{Results}

\subsection{Description of the problem}
\label{subsec:2.1}

Given an integer $M>1$, let $\Omega_M$ be the hierarchical group defined in the Introduction (see (\ref{eq:1.5}),(\ref{eq:1.6})). On $\Omega_M$ we consider a system of continuous time random walks
defined in the following way. At time $t=0$, at each $x\in \Omega_M$ there are $\theta_x$ particles, where $\{\theta_x\}_{x\in M}$ are independent copies of a random variable $\theta$ with values in $\{0,1,2,\ldots\}$, such that $E\theta^3<\infty$. A typical case  is that $\theta$ is a Poisson random variable. Such a system will be called Poisson. Each particle evolves according to a (continuous time) $c$-random walk independently of the others, with the same $c$ (see  (\ref{eq:1.7})). 
Sometimes we assume additionaly that 
the particles branch at rate $V$ independently, and the branching is critical and binary ($0$ or $2$ particles with equal probabilities). The evolution of the system is described by the empirical process 
$(N_t)_{t\geq 0}$, where 
$$N_t(A)=\,\hbox{\rm number of particles in the set}\, A\subset 
\Omega_M\,
\hbox{\rm at time}\, t.$$

For $T>0$, we define the occupation time fluctuation process
\begin{equation}
\label{eq:2.1}
X_T(t)=\frac{1}{F_T}\int^{Tt}_0(N_r-EN_r)dr,\quad t\geq 0,
\end{equation}
where $F_T$ is a suitable deterministic norming. We are interested in the asymptotics of the law of $X_T$ as $T\to\infty$. Unlike the system of  $\alpha$-stable motions in ${\mathbb R}^d$, in the most interesting cases a limit of $X_T$ 
exists only for special sequences $T_n\to\infty$.

\subsection{Convergence results}
\label{subsec:2.2}

We begin by defining the processes that will appear in the limits; their existence will follow from the corresponding convergence results and their properties will be given in the next subsection.

Denote
\begin{equation}
\label{eq:2.2}
a=\frac{c}{M},\quad b=\frac{M^2-c}{M(M-1)},
\end{equation}
(our $b$ is $ba$ in the notation of \cite{DGW1}).

{\it Oscillatory fractional Brownian motion (ofBm)} with parameter 
$H\in({1}/{2},1)$ is a centered Gaussian process 
$\xi^H=(\xi^H_t)_{t\geq 0}$ with covariance given by  (\ref{eq:1.3}) (hence it clearly has stationary increments), where
\begin{equation}
\label{eq:2.3}
h_t=\sum_{j\in{\mathbb Z}}(ba^jt)^{-2H}(e^{-ba^jt}-1+ba^jt),\quad t>0.
\end{equation}
Convergence of this series is easy to check, and
\begin{equation}
\label{eq:2.4}
0<C_1(a,H)<h_t<C_2(a,H)<\infty,
\end{equation}
for all $t>0$, with
\begin{equation}
\label{eq:2.4a}
C_1(a,H)=\frac{a^{2H-1}\Gamma(2-2H)}{(1-a^{2H-1})2H},\quad C_2(a,H)=a^{1-2H}C_1(a,H),
\end{equation}
(see \cite{DGW3}, (4.2.8) and Proposition 3.1.1). We additionally define $h_0=0$. Note that
\begin{equation}
\label{eq:2.5}
h_{at}=h_t,\quad t\geq 0,
\end{equation}
which justifies the qualifications ``oscillatory'' and ``periodic in logarithmic scale''. Observe that 
$h_t$ is infinite 
for $H\in (0,{1}/{2}]$. The constant $b$ in (\ref{eq:2.2}) appears in the limit, but it is irrelevant for the existence and properties of the process defined by (\ref{eq:1.3}) and (\ref{eq:2.3}); one may take any $b>0$.

{\it Oscillatory sub-fractional Brownian motion (osfBm)} with parameter 
$H\in ({1}/{2},1)$ is a centered Gaussian process $\zeta^H=(\zeta^H_t)_{t\geq 0}$ with covariance given by (\ref{eq:1.4}) and $h_t$ defined by (\ref{eq:2.3}).

{\it Oscillatory negative sub-fractional Brownian motion (onsfBm)} with parameter $H\in (1,{3}/{2})$ is centered Gaussian process $\zeta^H=(\zeta^H_t)_{t\geq 0}$ with covariance of the form (\ref{eq:1.4}) with $h$ replaced by
\begin{equation}
\label{eq:2.6}
\widetilde{h}_t=\sum_{j\in{\mathbb Z}}(ba^jt)^{-2H}\left(e^{-ba^jt}-1+ba^jt-
\frac{1}{2}(ba^jt)^2\right),\quad t>0.
\end{equation}
Using $e^{-x}-1+x-\frac{x^2}{2}=-\int^x_0\int^y_0\int^z_0e^{-v}dvdzdy,\, x>0$, it is easy to check that $\widetilde{h}_t$ is well defined and $|\widetilde{h}_t|$ satisfies (\ref{eq:2.4}) with
\begin{equation}
\label{2.6a}
C_1(a,H)=\frac{\Gamma(4-2H)}{2H(2H-1)(2H-2)(a^{2H-3}-1)},\,\, 
C_2(a,H)=a^{2H-3}C_1(a,H).
\end{equation}
$\widetilde{h}$ obviously satisfies (\ref{eq:2.5}) and $\widetilde{h}_t<0$. Note that putting formally $\widetilde{h}_t\equiv -1$ we obtain the negative sub-fractional Brownian motion considered in \cite{BGT7}.

\begin{remark}
\label{r:2.1}
{\rm It is useful to observe that the covariances of both {\it osfBm} and 
{\it onsfBm} can be written in a unified form:
\begin{equation}
\label{eq:2.7}
\int^t_0\int^s_0\int^{u\wedge v}_0\sum_{j\in{\mathbb Z}}e^{-ba^j(u+v-2r)}(ba^j)^{3-2H}drdudv.
\end{equation}}
After computing the integrals, \eqref{eq:2.7} takes the form
\begin{equation*}
 \sum_{j\in\Z} (A_j-B_j), \qquad A_j, B_j\ge 0, 
\end{equation*}
where the series $\sum_j{A_j}$, $\sum_j B_j$ converge if $H\in(1/2,1)$
(hence it can be written as ${\sum_jA_j-\sum_j B_j}$, yielding \eqref{eq:1.4} with $h$ of the form \eqref{eq:2.3})
and diverge if $H\in (1,3/2)$. In the latter case it is necessary to take $\sum_j(A_j-C_j)$, $\sum_j(B_j-C_j)$ for appropriate $C_j$, hence the different forms of $h$ and $\tilde h$.
For example, the existence of osfBm and onsfBm follows immediately from 
(\ref{eq:2.7}).
The expression (\ref{eq:2.7}) is finite and positive definite also for $H=1$, and the corresponding Gaussian process is an analogue of the process obtained in
(2.14) of  \cite{BGT7}, and it is discussed in \cite{BGT5} (formulas (3.5), (3.6)).
 However, in this case (\ref{eq:2.7}) cannot be given a form like (\ref{eq:1.4}).
\end{remark}

Another process that appears as an ingredient of the occupation time fluctuation limit is a centered Gaussian process $\eta^H$ with covariance
\begin{equation}
\label{eq:2.8}
E\eta^H_s\eta^H_t=(s+t)^{2H}h_{s+t}-t^{2H}h_t-s^{2H}h_s,\,\,s,t\geq 0,
\end{equation}
where $h$ is given by \eqref{eq:2.3} and $H\in(1/2,1)$. 
This is an oscillatory analogue of a process studied in \cite{LN,BB,RT} and, independently, in \cite{BGT9}.

We now give the convergence results. In what follows,  $K$ is a constant whose value can be calculated explicitely in each case.
This is our oscillatory analogue of a process introduced in \cite{LN} for $H\in (0,1/2)$, and extended in \cite{BB,RT} for $H\in (1/2,1)$.
\begin{theorem}
\label{t:2.2} Consider a system of $c$-random walks without branching and $c<1$ ($\gamma <0$, see \eqref{eq:1.9}), and a general initial condition determined 
by  $\theta$ (see Subsection \ref{subsec:2.1}).
Let $X_T$ be defined by \eqref{eq:2.1} with
\begin{equation}
\label{eq:2.9}
F_T=T^{{(1-\gamma)}/{2}}.
\end{equation}
Then, for 
\begin{equation}
\label{eq:2.10}
T_n =a^{-n},\quad n=1,2,\ldots,
\end{equation}
and any $\varphi\in {\cal B}_b$,
\begin{equation}
\label{eq:2.11}
\langle X_{T_n},\varphi\rangle
\Rightarrow_{\kern-.35cm _c}\,\,
K(\sqrt{2E\theta}
\zeta^H+\sqrt{{\rm Var}\theta}\eta^H)\sum_{x\in\Omega_M}\varphi(x),
\end{equation}
with $\zeta^H$ and $\eta^H$ independent, and 
$H={(1-\gamma})/{2}\in ({1}/{2},1)$ (i.e., $\zeta^H$ is osfBm).
\end{theorem}

\begin{remark}
\label{r:2.3}{\rm 
(a) If the distribution of $\theta$ is Poisson or, more generally, $E\theta={\rm Var}\theta$, then the limit process in (\ref{eq:2.11}) is $K_1(\varphi)\xi^H$. If $\theta$ is deterministic, then this limit is $K_2(\varphi)\zeta^H$.

\noindent
(b) This result is an analogue of Theorem 2.2 in \cite{BGT9}. It exhibits the relevant role of osfBm for the model.

\noindent
(c) From the proof it will be seen that instead of $T_n$ of the form (\ref{eq:2.10}) one can take any sequence $(T_n)_n$ such that $T_n a^{kn}$ has a positive limit for some $k=1,2,\ldots$. Then there exists a unique $m_0\in{\mathbb Z}$ and $\kappa\in [1,1/a)$ such that $T_n a^{kn}a^{m_0}\to \kappa$, and the limit of $X_{T_n}$ has the form (\ref{eq:2.11}), where in the definition of $h$ (see (\ref{eq:2.3})) $b$ is replaced by $\kappa b$.

\noindent
(d) We have formulated the convergence results for a fixed test function $\varphi\in {\cal B}_b$ only, but as $\Omega_M$ is countable and discrete, it is easily seen that the (weak) convergences take place also in the space of continuous processes with values in the dual ${\cal B}_b'$.
\medskip

Remarks (c) and (d) apply also to the next two theorems.

\medskip
We now pass to the branching case.}
\end{remark}

\begin{theorem}
\label{t:2.4}
Consider a Poisson system of $c$-random walks with binary branching at rate $V$. We assume that $1<c<\sqrt{M}\,\, (0<\gamma<1)$. Then, with
\begin{equation}
\label{eq:2.12}
F_T=T^{{(2-\gamma})/{2}},
\end{equation}
$T_n$ of the form \eqref{eq:2.10} and any $\varphi\in {\cal B}_b$,
\begin{equation}
\label{eq:2.13}
\langle X_{T_n},\varphi\rangle\Rightarrow_{\kern-.35cm _c}\,\,
K\zeta^H\sum_{x\in\Omega_M}\varphi(x),
\end{equation}
where $\zeta^H$ is osfBm with $H={(2-\gamma})/{2}\in ({1}/{2},1)$.
\end{theorem}

\begin{remark}
\label{r:2.5}{\rm 
This theorem is an analogue of Theorem 2.2 for $\alpha$-stable motions in \cite{BGT1}. The latter result was extended to general initial conditions of homogeneous type with the same limit process in \cite{BGT9}. However this extension was not immediate. We are sure that (\ref{eq:2.13}) also holds for any initial condition described in Subsection 2.1, but we have not made detailed calculations.}
\end{remark}

In order to obtain onsfBm it is necessary to add high density to the previous model.

\begin{theorem}
\label{t:2.6}
Consider a Poisson system of $c$-random walks with binary branching at rate $V$. We assume that $c<1\,\,(\gamma<0)$ and the Poisson intensity of the initial distribution is $H_T>0$ such that
\begin{equation}
\label{eq:2.14}
\lim_{T\to\infty}T^\gamma H_T=\infty.
\end{equation}
Then, with
\begin{equation}
\label{eq:2.15}
F_T=T^{{(2-\gamma})/{2}} H^{{1}/{2}}_T,
\end{equation}
$T_n$ of the form \eqref{eq:2.10} and any $\varphi\in {\cal B}_b$,
\begin{equation}
\label{eq:2.16}
\langle X_{T_n},\varphi\rangle 
\Rightarrow_{\kern-.35cm _c}\,\,
K\zeta^H\sum_{x\in\Omega_M}\varphi(x),
\end{equation}
where $\zeta^H$ is onsfBm with $H={(2-\gamma})/{2}\in (1,{3}/{2})$.
\end{theorem}

\begin{remark}
\label{r:2.7}{\rm 
This theorem is an analogue of Theorem 2.2(a) of \cite{BGT7} 
(see also Remark 2.3(a) and Proposition 2.5 therein). In \cite{BGT7} a convergence result was also given for $\alpha=d$, which in our case corresponds to $c=1\,\, (\gamma=0$, see 
(\ref{eq:1.8}), (\ref{eq:1.9})). For such $c$ one can obtain (\ref{eq:2.16})
as well with $\zeta^H$ whose covariance is of the form (\ref{eq:2.7}) with $H=1$ (see Remark 2.1).}
\end{remark}

We now pass to the case of ``large'' $\gamma$.

\begin{proposition}
\label{p:2.8} 
Consider a Poisson system of $c$-random walks without branching with \\
$c>1\,\,(\gamma>0)$. Let $F_T=\sqrt{T}$. Then
$$X_T
\Rightarrow_{\kern-.35cm _f}\,\,
X\quad{\rm as}\quad T\to\infty,$$
where $X$ is a centered Gaussian process with values in the dual ${\cal B}_b'$ with covariance
\begin{equation}
\label{eq:2.17}
E\langle X(s),\varphi\rangle\langle X(t),\psi\rangle=2E\theta(s\wedge t)\langle G\varphi,\psi\rangle,
\end{equation}
where $G$ is the potential operator of the particle motion process.
\end{proposition}

\begin{remark}
\label{r:2.9} 
{\rm
(a) We formulate explicitly the result for the Poisson model without branching only, since in other cases the limit processes are similar. It seems clear that analogously as for the $\alpha$-stable motions in $\mathbb R^d$ (see \cite{BGT9}, complete version), replacing Poisson random measure by an arbitrary initial 
distribution described in Subsection 2.1 does not change the limit. In the branching case with $\gamma >1\,\, (c>\sqrt{M})$, a result of \cite{DGW1} (Theorem 
2.2.2), and \cite{BGT2} and \cite{BGT9} (complete version, $\alpha$-stable motions) indicate that, again with $F_T=\sqrt{T}$, the limit process is the sum of two independent centered Gaussian processes, one of the form (\ref{eq:2.17}), and the other with covariance 
$E\theta\, V(s\wedge t)\langle G^2\varphi,\psi\rangle$. In the critical cases, i.e., $\gamma=0$ (no branching) and $\gamma=1$ (branching), with $F_T=\sqrt{T\log T}$, the limits are of the form $K\mu\beta$, where $\beta$ is the standard one dimensional Brownian motion, and $\mu$ is the counting measure on $\Omega_M$. As these results do not seem particularly interesting, 
we dispense with the calculations.

\noindent
(b) So we see that for ``large'' $\gamma$, analogously as for $\alpha$-stable motions, the limit processes have independent stationary increments. Moreover, the form of the covariances is expressed in terms of the potential operator $G$ in the same way. However, an interesting feature of the process with covariance 
(\ref{eq:2.17}) is that its spatial structure is completely different from the corresponding limit for $\alpha$-stable motions. We will see this in 
Theorem \ref{t:2.19} and Remark 2.22.

\noindent
(c) Observe that for ``large'' $\gamma$ no oscillations occur, therefore to obtain convergence there is no need to take special sequences $T_n$ (see (\ref{eq:2.10})).}
\end{remark}

\subsection{Properties of the limit processes}
\label{subsec:2.3}

\subsubsection{Different forms of the covariances }
\label{subsub:2.3.1}

\begin{proposition}
\label{p:2.10}

(a) Let $h$ be given by \eqref{eq:2.3} and $H\in({1}/{2},1)$, then
\begin{equation}
\label{eq:2.18}
t^{2H}h_t=\int^\infty_0\frac{\sin^2(tu)}{u^2}2\sigma(2u)du,
\end{equation}
where
\begin{equation}
\label{eq:2.19}
\sigma(u)=\frac{1}{\pi}\sum_{j\in{\mathbb Z}}
\frac{(ba^j)^{3-2H}}{(ba^j)^2+u^2}.\end{equation}
(b) For ofBm $\xi^H$ we have
\begin{equation}
\label{eq:2.20}
E\xi^H_s\xi^H_t=\int^\infty_0\frac{(1-\cos(su))(1-\cos(tu))+\sin(su)\sin(tu)}{u^2}\sigma(u)du.
\end{equation}
(c) For osfBm and onsfBm 
$\zeta^H, H\in ({1}/{2},1)\cup (1,{3}/{2})$, we have
\begin{equation}
\label{eq:2.21}
E\zeta^H_s\zeta^H_t=\int^\infty_0\frac{(1-\cos(su))(1-\cos(tu))}{u^2}\sigma(u)du.
\end{equation}
\end{proposition}

\begin{remark}
\label{r:2.11}{\rm 
(a) Formulas (\ref{eq:2.20}) and (\ref{eq:2.21}) for $H\in ({1}/{2},1)$ follow immediately from (\ref{eq:2.18}), but the fact that (\ref{eq:2.21}) also holds for onsfBm $(H\in(1,{3}/{2}))$ is not so immediate since the covariance is expressed in terms of $\widetilde{h}$ given by (\ref{eq:2.6}), and not $h$.
The covariance (\ref{eq:1.3}) with the representation 
\eqref{eq:2.18}-\eqref{eq:2.19} (extending $t^{2H}h_t$ by symmetry to $t<0$) is of the form related to ``helical arcs'' studied in \cite{NS} (see also \cite{AL,Kah}).

\noindent
(b) From this proposition it is clear that ofBm, osfBm and onsfBm can be defined for all $t\in {\mathbb R}$, if in (\ref{eq:1.3}) and (\ref{eq:1.4}) we put $|t|, |s|, |t+s|$. Then the process
$$\left(\frac{\xi^H_t+\xi^H_{-t}}{\sqrt{2}}\right)_{t\in{\mathbb R}}$$
is an osfBm (this relationship for fBm and sfBm was observed in \cite{BGT0}).}
\end{remark}

The next simple proposition exhibits another aspect of analogies between the structures of oscillatory processes and their non-oscillatory counterparts. 
The corresponding formulas for the non-oscillatory processes turned out to be useful in \cite{BT}.

\begin{proposition}
\label{p:2.12}

(a) For ofBm $\xi^H$ we have
\begin{equation}
\label{eq:2.22}
E\xi^H_s\xi^H_t=\int_{\mathbb R^2}\UNO_{[0,s]}(x)\UNO_{[0,t]}(y)\sum_{j
\in{\mathbb Z}}(ba^j)^{2-2H}e^{-|x-y|ba^j}dxdy.
\end{equation}
(b) For osfBm $\zeta^H,H\in ({1}/{2},1)$, we have
\begin{equation}
\label{eq:2.23}
E\zeta^H_s\zeta^H_t=\frac{1}{4}\int_{\mathbb R^2}(\UNO_{[0,s]}-\UNO_{[-s,0]})(x)(\UNO_{[0,t]}-\UNO_{[-t,0]})(y)\sum_{j\in{\mathbb Z}}(ba^j)^{2-2H}e^{-|x-y|ba^j}dxdy.
\end{equation}
(c) For onsfBm $\zeta^h, H\in(1,{3}/{2})$, we have
\begin{equation}
\label{eq:2.24}
E\zeta^H_s\zeta^H_t=\frac{1}{4}\sum_{j\in {\mathbb Z}}\int_{\mathbb R^2}(\UNO_{[0,s]}-\UNO_{[-s,0]})(x)(\UNO_{[0,t]}-\UNO_{[-t,0]})(y)(ba^j)^{2-2H}e^{-|x-y|ba^j}dxdy.
\end{equation}
\end{proposition}

\begin{remark}
\label{l:2.13}{\rm 
(a) (\ref{eq:2.24}) holds also for osfBm, but for onfsBm the sum cannot be put under the integral.

\noindent
(b) For fractional Brownian motion and sub-fractional Brownian motion with 
$H\in({1}/{2},1)$ there are analogues of (\ref{eq:2.22}) and (\ref{eq:2.23}) with kernel $C|x-y|^{2H-2}$ in \cite{BT} (Corollary 2.12 and Remark 2.13(a)).}
\end{remark}

\subsubsection{Representations in terms of a sequence of independent Brownian motions}
\label{subsub:2.3.2}

\begin{proposition}
\label{p:2.14}
Let $\{\beta^j\}_{j\in\mathbb Z}$ be independent standard Brownian motions and $\{\nu_j\}_{j\in\mathbb Z}$ be an i.i.d. family of standard normal random variables, independent of $\{\beta^j\}_{j\in\mathbb Z}$.

\noindent
(a) Let $H\in({1}/{2},1)$. The process
\begin{equation}
\label{eq:2.25}
\sum_{j\in\mathbb Z}(ba^j)^{-H}(1-e^{-ba^jt})
\nu_j+\sqrt{2}\sum_{j\in\mathbb Z}(ba^{j})^{{1}/{2}-H}
\int^t_0(1-e^{-ba^j(t-r)})d\beta^j_r,\quad t\geq 0,
\end{equation}
is an ofBm.

\noindent
(b) The process
\begin{equation}
\label{eq:2.26}
\sum_{j\in\mathbb Z}(ba^j)^{{1}/{2}-H}\int^t_0(1-e^{-ba^j(t-r)})d\beta^j_r,\quad t\geq 0,
\end{equation}
is an osfBm if $H\in({1}/{2},1)$ and an onsfBm if 
$H\in (1,{3}/{2})$.
\end{proposition}

\begin{remark}
\label{r:2.15}{\rm 
(a) the first series in (\ref{eq:2.25}) is the process $\eta^H$ defined by
(\ref{eq:2.8}), and (\ref{eq:2.25}) corresponds to the representation of ofBm resulting from (\ref{eq:2.11}) for $E\theta={\rm Var}\theta$ (see Remark 2.3(a)).

\noindent
(b) It is not clear how to represent ofBm and osfBm in terms of just one 
stochastic integral. For fBm there are representations with one stochastic integral, see \cite{Nu} and references therein; for sfBm, see \cite{BGT0}.}
\end{remark}

\subsubsection{Integral representation of onsfBm}
\label{subsub:2.3.3}

Negative sub-fractional Brownian motion with parameter $H$ can be represented as $C\!\int^t_0\varrho_u du$ where $\varrho$ is the odd part in the decomposition of fBm with parameter $H-1$ considered in \cite{DZ} (see \cite{BGT5}, Theorem 3.1). A similar fact is true for onsfBm with an oscillatory 
analogue of $\varrho$.

\begin{proposition}
\label{p:2.16}

OnsfBm $(H\in (1,{3}/{2}))$ has a representation
\begin{equation}
\label{eq:2.27}
\zeta^H_t=\int^t_0\varrho^{H-1}_udu,\quad t\geq 0,
\end{equation}
where $\varrho^U$ is a continuous centered Gaussian process with covariance
\begin{equation}
\label{eq:2.28}
E\varrho^U_s\varrho^U_t=\frac{1}{2}\left[(s+t)^{2U}
\tilde{\tilde h}^U_{s+t}-|s-t|^{2U}
\tilde{\tilde h}^U_{|s-t|}
\right] ,\,\, U\in (0,{1}/{2}),
\end{equation}
and
\begin{equation}
\label{eq:2.29}
\tilde{\tilde h}^U_t=\sum_{j\in\mathbb Z}(ba^jt)^{-2U}(1-e^{-ba^jt}),\quad t>0.
\end{equation}
\end{proposition}

The integrand process in  the analogue of (\ref{eq:2.27}) for the 
nsfBm has covariance of the form (\ref{eq:2.28}) with 
$\tilde{\tilde h}\equiv 1$.

\subsubsection{Bounds for second moments of increments and H\"older continuity }
\label{subsub:2.3.4}

\begin{proposition}
\label{p:2.17}

Let $\vartheta^H$ denote either one of the processes $\xi^H$ or $\zeta^H$, then 
\begin{equation}
\label{eq:2.30}
C_1|t-s|^{(2H)\wedge 2}\leq E(\vartheta^H_t-\vartheta^H_s)^2\leq C_2|t-s|^{(2H)\wedge 2},
\end{equation}
where $C_1$ and $C_2$ are positive constants that depend on $\xi^H$ or 
$\zeta^H$, respectively; in particular for osfBm 
$\xi^H$ these constants are given by \eqref{eq:2.4a}.
\end{proposition}
\vglue.5cm
Hence ofBm and osfBm have versions whose paths are locally H\"older continuous with any exponent smaller than $H$. On the other hand, Proposition 2.16 implies that the paths of onsfBm are continuously differentiable.

\subsubsection{Semimartingale/non-semimartingale property}
\label{subsub:2.3.5}

From (\ref{eq:2.30}) and Lemma 2.1 in \cite{BGT} it follows that ofBm and osfBm are not semimartingales. It can be also shown that the centered Gaussian process with covariance (\ref{eq:2.7}) for $H=1$ is not semimartingale either. On the other hand, onsfBm is obviously a semimartingale by (\ref{eq:2.27}).
The process $\eta^H$ defined by \eqref{eq:2.8} is also a semimartingale, since from \eqref{eq:2.25} and Remark \ref{r:2.15} (a) it follows that it can be written in a form similar to \eqref{eq:2.27}.

\subsubsection{Weighted quadratic variation}
\label{subsub:2.3.5a}

The result of this subsection a first attempt to clarify how the oscillatory character of the covariance function is reflected in the behavior of paths. From \eqref{eq:2.30} it follows that the quadratic variation of the paths of ofBm and osfBm is $0$ on each interval $[0,\tau]$.  We now study the rate of convergence 
to $0$. An analogous problem  for fBm and sfBm was solved in \cite{Gl,N,T1}.

\begin{proposition}
 \label{p:2.17a}
 Fix $T>0$, $H\in(1/2, 1)$ and let $\vartheta$ be either ofBm or osfBm

\noindent
(a) For $\varepsilon>0$, let 
\begin{equation}
 \label{eq:2.30a}
 V_\varepsilon(\vartheta)=\varepsilon^{-2H}\int_0^T(\vartheta_{t+\varepsilon}-\vartheta_t)^2dt.
\end{equation}
Then
\begin{equation}
 \label{eq:2.30b}
 \lim_{\varepsilon\to 0}\frac 1{h_\varepsilon}V_\varepsilon(\vartheta)=T
\end{equation}
in $L^2$, where $h_\varepsilon$ is given by \eqref{eq:2.3}. Moreover, if $(\varepsilon_n)_{n=1,2,\ldots}$ is such that 
\begin{equation}
 \label{eq:2.30c}
 \sum_{n=1}^\infty \varepsilon_n^\delta<\infty\qquad \textrm{for some }\quad 0<\delta<(4-4H)\wedge 1,
\end{equation}
then for $\varepsilon=\varepsilon_n$ the convergence in 
\eqref{eq:2.30b} as $n\to\infty$ is almost sure as well.

\noindent
(b) The same assertion holds for 
\begin{equation}
U_\varepsilon(\vartheta)=\varepsilon^{1-2H}
\sum_{k=1}^{\lfloor  T/\varepsilon\rfloor} (\vartheta_{k\varepsilon}-\vartheta_{(k-1)\varepsilon})^2,
\end{equation}
(where $\lfloor\cdot\rfloor$ stands for integer part).
\end{proposition}

\begin{remark}
 \label{r:2.17b}
 (a) For a sequence $\varepsilon_n=\kappa a^{n}$, $\kappa\in[1, 1/a)$, \eqref{eq:2.30c} is obviously satisfied, and then we have
 \begin{equation*}
  \lim_{n\to\infty}V_{\varepsilon_n}(\vartheta)= \lim_{n\to\infty}U_{\varepsilon_n}(\vartheta)=T h_{\kappa}
 \end{equation*}
in $L^2$ and almost surely. This result exhibits the oscillatory (in a sense) character of the paths.

\noindent
(b) In \cite{T1} the $p$-variation of sfBm is studied. It seems that similar methods would allow to investigate ofBm and osfBm.
\end{remark}

\subsubsection{Non-Markov}
\label{subsub:2.3.6}

None of the covariances of the oscillatory processes has the triangular property \cite{Kal} (Proposition 13.7), hence those processes are not Markov.

\subsubsection{Self-similarity}
\label{subsub:2,3,7}

$\xi^H$ and $\zeta^H$ are not self-similar, but we have
$$(\xi^H_{\lambda t})_{t\geq 0}\stackrel{d}{=} (\lambda^H\xi^H_t)_{t\geq 0}$$
for $\lambda=a^j, j\in{\mathbb Z}$, by (\ref{eq:2.5}) and (\ref{eq:1.3}). Moreover, for general $\lambda$,
if we write it as $\lambda=a^j\kappa,\\
 \kappa \in [1,{1}/{a}), j\in{\mathbb Z}$, then
$$(\xi^H_{\lambda t})_{t\geq 0}\stackrel{d}{=}(\lambda^H\xi^{\kappa,H}_t)_{t\geq 0},$$
where $\xi^{\kappa,H}$ is ofBm with $b$ replaced by $b\kappa$ (see (\ref{eq:2.3})). The same is true for $\zeta^H$.

\subsubsection{Long range dependence}
\label{eq:2.3.8}

\begin{proposition}
\label{p:2.18}

Let $0\leq u<v\leq s<t$, then

\noindent
(a)
\begin{eqnarray*}
0&<&C^H_1(u,v,s,t)\leq\liminf_{\tau\to\infty}\tau^{2-2H}E(\xi^H_v-\xi^H_u)(\xi^H_{t+\tau}-\xi^H_{s+\tau})\\
&\leq&\limsup_{\tau\to\infty}\tau^{2-2H}E(\xi^H_v-\xi^H_u)(\xi^H_{t+\tau}-\xi^H_{s+\tau})\leq C^H_2(u,v,s,t)<\infty,
\end{eqnarray*}
for  $H\in(1/2, 1)$,

\noindent
(b)
\begin{eqnarray*}
0&<&C^H_3(u,v,s,t)\leq\liminf_{\tau\to\infty}\tau^{3-2H}E(\zeta^H_v-\zeta^H_u)
\zeta^H_{t+\tau}-\zeta^H_{s+\tau})\\
&\leq&\limsup_{\tau\to\infty}\tau^{3-2H}E(\zeta^H_v-\zeta^H_u)
(\zeta^H_{t+\tau}-\tau^H_{s+\tau})\leq C^H_4(u,v,s,t)<\infty
\end{eqnarray*}
for $H\in({1}/{2},1)\cup(1,{3}/{2})$, where the $C^H_i(u,v,s,t)$ are positive constants for each $u,v,s,t$.
\end{proposition}

The orders of decay of the covariances of increments for the oscillatory processes are the same as the orders for their non-oscillatory counterparts (see \cite{BGT0,BGT5}, the result is well known for fBm).

\subsubsection{Asymptotics for the boundary values of $H$}
\label{subsub:2.3.9}

Both fBm and sfBm  are the standard Brownian motion for $H={1}/{2}$, and 
there is a continuous dependence on $H$. We have already mentioned that ofBm and osfBm are not defined for $H={1}/{2}$. Nevertheless, from 
\eqref{eq:2.4}, \eqref{eq:2.4a} and \eqref{eq:1.3} we have that
$$\sqrt{1-a^{2H-1}}\xi^H
\Rightarrow_{\kern-.35cm _f} 
\,\,K\beta \quad{\rm as}\quad H\searrow{1}/{2},$$
where $\beta$ is  standard Brownian motion. The same holds for osfBm $\zeta^H$.

Similarly, for onsfBm, from \eqref{eq:2.6}, \eqref{2.6a} and \eqref{eq:1.4} it follows that if 
$H\nearrow {3}/{2}$, then
$\sqrt{a^{2H-3}-1}\zeta^H$ converges in the sense of finite dimensional distributions to nsfBm with $H={3}/{2}$ (see 
\cite{BGT5,BGT7}).

\subsubsection{Spatial structure  for ``large'' $\gamma$}
\label{subsub:2.3.10}

We restrict our considerations here to a Poisson non-branching system with $\gamma>0$ ($c>1$). In the remaining cases the results will be either identical or similar (see Remark 2.9(a)). With no loss of generality we assume the Poisson parameter $E\theta=1$.

Let $B_n$ denote the ball in $\Omega_M$ with center $0$ and radius $n=0,1,\ldots$. Denote
\begin{equation}
\label{eq:2.31}
D=\frac{M-1}{Mb} \frac{c}{c-1},
\end{equation}
where $b$ is given by (\ref{eq:2.2}).

\begin{theorem}
\label{t:2.19}
Let $X$ be the limit process given by  Proposition \ref{p:2.8} and denote
\begin{equation}
\label{eq:2.32}
\widetilde{\vartheta}_t=\langle X(1),\UNO_{B_t}\rangle,\,\, t=0,1,\ldots .
\end{equation}
Then $\widetilde{\vartheta}$ is a centered Gaussian process with covariance
\begin{equation}
\label{eq:2.33}
E\widetilde{\vartheta}_s\widetilde{\vartheta}_t=2D M^{s\wedge t}
\left(\frac{M}{c}\right)^{s\vee t}.
\end{equation}
\end{theorem}

\begin{remark}
\label{r:2.20}
{\rm 
(a) From (\ref{eq:2.33}) it follows that the process $\widetilde{\vartheta}$ and its continuous time interpolation resulting from (\ref{eq:2.33}) are Markov 
(\cite{Kal}, Proposition 13.7). This is in sharp contrast with the $\alpha$-stable case. Recall that  $\gamma>0$ corresponds to $1=d>\alpha$ for the Poisson system of $\alpha$-stable motions in $\mathbb R$ without branching. If $X$ is the occupation time fluctuation limit process (see \cite{BGT2} or \cite{DGW1}), then the counterpart of (\ref{eq:2.32}) is $\langle X(1),\UNO_{[-t,t]}\rangle$, and this process is the odd part of the fBm in the sense of \cite{DZ}, i.e., a centered Gaussian process with covariance of the form
$$K(|t+s|^{2H}-|t-s|^{2H}),\quad H=\frac{1+\alpha}{2}.$$
This process is not Markov. In the $\alpha$-stable case it makes sense to consider also the processes $\langle X(1), \UNO_{[0,t]}\rangle$ and $\langle X(1),
\UNO_{[0,t]}-\UNO_{[-t,0]}\rangle$, which are fBm and sfBm, respectively \cite{BT} (Corollary 2.12), see also \cite{LX}. So the results of \cite{BGT2} and \cite{BT} show that in this model, as $\gamma$ changes from ``small'' to ``large'' values, there is a phase transition in the sense that a fractional temporal structure and a uniform spatial structure changes into a simple temporal structure and a fractional spatial structure. In the hierarchical case we have a passage from oscillatory fractional temporal structure and uniform spatial structure to a simple temporal structure and a Markov (non-fractional) spatial structure.

\noindent
(b) It is easy to see that the process $\widetilde{\vartheta}$ can be obtained as a limit of $(\langle X_T(1),\UNO_{B_t}\rangle)_{t=0,1,\ldots}$.

\noindent
(c) The hierarchical structure of the state space $\Omega_M$ is the reason that the logarithmic scale is more natural. So, we consider the process $\vartheta$ defined by
\begin{equation}
\label{eq:2.34}
\vartheta_t=\widetilde{\vartheta}_{\log_{{1}/{a}}t},\,\, t>0.
\end{equation}
This is a centered Gaussian process, and by (\ref{eq:2.33})
its covariance is given by
\begin{equation}
\label{eq:2.35}
E\vartheta_s\vartheta_t=2D(s\wedge t)^{1+\gamma}(s\vee t).
\end{equation}
By (\ref{eq:2.34}), $\vartheta_t$ is defined for $t=a^{-n},n=0,1,2,\ldots$, but a Gaussian process with covariance (\ref{eq:2.35}) is well defined for any $t\geq 0$. The interpolation is also denoted by $\vartheta$. The usefulness of this type of logaritmic scaling and the relevance of the sequence $a^{-n}$ (see  
(\ref{eq:2.10})) are also seen in \cite{BGT10}.
\vglue.5cm
We now collect some properties of the process $\vartheta$ defined by \eqref{eq:2.34}, which follow immediately from \eqref{eq:2.35}.}
\end{remark}

\begin{proposition}
\label{p:2.21}

The centered Gaussian process $\vartheta$ with covariance \eqref{eq:2.35}, with $\gamma>0$, has the following properties:

\noindent
(i) it is Markov,

\noindent
(ii) it is self-similar with index $2+\gamma$,

\noindent
(iii) it can be represented as
$$\vartheta_t=\sqrt{2D\gamma} t\int^t_0u^{({\gamma-1})/{2}}d\beta_u,$$
where $\beta$  is a standard Brownian motion,

\noindent
(iv) it is of diffusion type, satisfying the  stochastic differential equation
$$d\vartheta_t=\frac{1}{t}\vartheta_tdt+\sqrt{2D\gamma}
t^{({\gamma+1})/{2}}d\beta_t.$$
\end{proposition}

\section{Proofs}
We start by recalling some facts on $c$-random walks. First observe that \eqref{eq:1.9} and \eqref{eq:2.2} imply

\begin{equation}
\label{eq:3.1}
\frac{1}{M}=a^{\gamma+1}.
\end{equation}
This permits to write the formula for the transition probability, given, e.g., 
in (3.1.3) of  \cite{DGW1}, in the following form
\begin{equation}
p_t(x,x+y)=p_t(0,y)=(\UNO_{\{0\}}(y)-1)a^{|y|(\gamma+1)}e^{-ba^{|y|-1}t}
+\frac{M-1}{M}\sum_{j=|y|}^\infty a^{j(\gamma+1)}e^{-ba^jt}\label{eq:3.2}
\end{equation}
(see \eqref{eq:2.2} for our definition of $b$).
This and an easy estimate,
\begin{equation}
\label{eq:3.3}
0<C_1\le\sum_{j\in\Z}(ba^jt)^{\gamma+1}e^{-ba^jt}<C_2,\,\, \gamma>-1,\,\, t>0,
\end{equation}
yield
\begin{equation}
\label{eq:3.4}
p_t(0,y)\le p_t(0,0) \le C_3(1\wedge t^{-(\gamma+1)}).
\end{equation}
Denote
\begin{equation}
\label{eq:3.5}
G_t\varphi(x)=\int_0^t{\cal T}_s\varphi(x)ds, \quad \varphi\in {\cal B}_b,
\end{equation}
\begin{numcases}{G_t\varphi(x)\le}
C(\varphi)t^{-\gamma}& if \ $-1<\gamma<0$
\label{eq:3.5a}\\
C(\varphi)\log t& if \ $\gamma=0$
\label{eq:3.5b}\\
C(\varphi)& if \ $\gamma>0.$
\label{eq:3.5c}
\end{numcases}
%
%
\vglue.5cm
\noindent
{\bf Proof of Theorem \ref{t:2.2}}
\vglue.4cm
\noindent
Recall that $\theta_x$ is the initial number of particles at $x\in\Omega_M$, and $\{\theta_x\}_{x\in\Omega_M}$ are independent copies of $\theta$. Denote
\begin{equation}
\label{eq:3.6}
p_k=P(\theta=k),\quad k=0,1,\dots
\end{equation}
Let $N^{(x)}$ be the empirical process of the system starting from $\theta_x$ particles at $x$. By assumption, the 
$N^{(x)}$, $x\in\Omega_M$, are independent and
\begin{equation}
\label{eq:3.7}
N=\sum_{x\in\Omega_M}N^{(x)}.
\end{equation}
This representation permits to use the central limit theorem in the proof of convergence of finite dimensional distributions.

For simplicity, we will prove the theorem for $\varphi=\UNO_{\{0\}}$. It will be clear that the same argument can be carried out in the general case.

Let $\{z^{x,j}\}_{x\in\Omega_M,j=1,2,\dots}$ be independent $c$-random walks. Then
\begin{equation}
\label{eq:3.8}
\langle N^{(x)}_t,\varphi\rangle=\sum_{j=1}^{\theta_x}\varphi(x+z^{x,j}_t).
\end{equation}

First we show convergence of covariances. Fix $t\ge s\ge 0$. We have 

\begin{equation}
\label{eq:3.9}
{\rm Cov}(\langle N_s^{(x)},\varphi\rangle, \langle N_t^{(x)},\varphi\rangle)=I-I\!\!I,
\end{equation}
where 
\begin{eqnarray}
I&=&E\langle N_s^{(x)},\varphi\rangle \langle N_t^{(x)},\varphi\rangle \nonumber\\
&=&\sum_{k=1}^\infty p_k
\sum_{\substack{i,j=1 \\ i\ne j}}^k
E\varphi(x+z_s^{x,j})E\varphi(x+z_t^{x,i})+\sum_{k=1}^\infty p_k\sum_{j=1}^k E(\varphi(x+z_s^{x,j})\varphi(x+z_t^{x,j}))\nonumber\\
&= &(E\theta^2-E\theta)p_s(0,x)p_t(0,x)+E\theta\, p_s(0,x)p_{t-s}(0,0),\label{eq:3.10}
\end{eqnarray}
and
\begin{equation}
\label{eq:3.11}
I\!\!I=E\langle N_s^{(x)},\varphi\rangle E\langle N_t^{(x)},\varphi\rangle=(E\theta)^2\, p_t(0,x)p_s(0,x).
\end{equation}
This, \eqref{eq:3.7} and the Chapman-Kolmogorov equation imply
$$
{\rm Cov}(\langle N_s,\varphi\rangle, \langle N_t,\varphi\rangle)=E\theta\, p_{t-s}(0,0)+({\rm Var}\, \theta-E\theta)p_{t+s}(0,0).
$$
Hence, by \eqref{eq:2.1}, 

\begin{equation}
\label{eq:3.12}
{\rm Cov}(\langle X_{T_n}(s),\varphi\rangle,\langle X_{T_n}(t),\varphi\rangle)=E\theta A_n+({\rm Var}\, \theta-E\theta)B_n,
\end{equation}
where

\begin{equation}
\label{eq:3.13}
A_n=\frac{T_n^2}{F_{T_n}^2}\int_0^s\int_0^tp_{T_n|u-v|}(0,0)dudv,
\end{equation}

\begin{equation}
\label{eq:3.14}
B_n=\frac{T_n^2}{F_{T_n}^2}\int_0^s\int_0^tp_{T_n(u+v)}(0,0)dudv.
\end{equation}

From \eqref{eq:2.9}, \eqref{eq:2.10} and \eqref{eq:3.2} it follows that

$$
A_n=\frac{M-1}{M}\int_0^s\int_0^t\sum_{j=0}^\infty a^{(j-n)(\gamma+1)}e^{-ba^{j-n}|u-v|}dudv.
$$

Using \eqref{eq:3.3} and $\gamma<0$ we have

\begin{equation}
\label{eq:3.15}
\lim_{n\to\infty}A_n=\frac{M-1}{M}\int_0^s\int_0^t\sum_{j\in\zetita} a^{j(\gamma+1)}e^{-ba^j|u-v|}dudv<\infty.
\end{equation}
Analogously,

\begin{equation}
\label{eq:3.16}
\lim_{n\to\infty}B_n=\frac{M-1}{M}\int_0^s\int_0^t\sum_{j\in\zetita} a^{j(\gamma+1)}e^{-ba^j(u+v)}dudv<\infty.
\end{equation}

Computing the integrals in \eqref{eq:3.15} and \eqref{eq:3.16}, by \eqref{eq:3.12} we obtain

\begin{equation}
\label{eq:3.17}
\lim_{n\to\infty}{\rm Cov}(\langle X_{T_n}(s),\varphi\rangle,\langle X_{T_n}(t),\varphi\rangle)=K^2(2E\theta E\zeta_s^H\zeta_t^H+ {\rm Var}\,\theta E\eta_s^H\eta_t^H),
\end{equation}
with $H=({1-\gamma})/{2}$ and $K^2=({M-1})/(Mb^{1+\gamma})$ (see \eqref{eq:1.4}, \eqref{eq:2.3} and \eqref{eq:2.8}).

Observe that existence of $\zeta^H$, $\eta^H$ and $\xi^H$ for 
$H\in(1/2,1)$ follows from \eqref{eq:3.17}. Namely, taking $\theta$ deterministic we obtain positive definiteness of \eqref{eq:1.4}. Next, taking $E\theta$ small and keeping ${\rm Var}\,\theta=1$ we see that \eqref{eq:2.8} is also positive definite. Finally, existence of $\xi^H$ is seen by taking $E\theta={\rm Var}\,\theta$.

Now we pass to the proof of convergence of finite dimensional distributions.

For fixed $a_1,\dots,a_m\in{\mathbb R}$
and $t_1,\dots,t_m\ge 0$, we write the sum $\sum_{k=1}^ma_k\langle X_{T_n}(t_k),\varphi\rangle$ with the help of \eqref{eq:3.7} and we apply the Lyapunov criterion. It is not hard to see, using $E\theta^3<\infty$, that it suffices to prove that

\begin{equation}
\label{eq:3.18}
S_{T_n}:= \frac{1}{F_{T_n}^3}\sum_{x\in\Omega_M}\left(\int_0^{T_n}\varphi(x+z_u^{x,1})du\right)^3\to 0\quad {\rm as}\ n\to\infty
\end{equation}
(see the end of Subsection 3.2 in \cite{BGT9}, or Subsection 3.3 in the complete version).

We have, using \eqref{eq:3.5} and \eqref{eq:1.10},
\begin{eqnarray}
S_T&=&\frac{6}{F_T^3}\sum_{x\in\Omega_M}\int_0^T\int_r^T\int_s^T{\cal T}_r(\varphi{\cal T}_{s-r}(\varphi{\cal T}_{t-s}\varphi))(x)dtdsdr\label{eq:3.18a}\nonumber\\
&\le&
\frac{6}{T^{{3(1-\gamma)}/{2}}}
\sum_{x\in\Omega_M}\int_0^T\varphi(x)G_T(\varphi G_T\varphi)(x)dr\nonumber\\
&=&\frac{6}{T^{({1}-{3\gamma)}/{2}}}G_T(\varphi G_T\varphi)(0)
\end{eqnarray}
(recall that $\varphi=\UNO_{\{0\}}$). By \eqref{eq:3.5a} we obtain
$$
S_T\le \frac{C}{T^{{(1+\gamma)}/{2}}},
$$
hence \eqref{eq:3.18} follows.

It remains to prove tightness of $\langle X_{T_n},\varphi\rangle$ in $C([0,\tau])$ for any $\tau>0$. From \eqref{eq:3.12}--\eqref{eq:3.14} it follows that,
 for $s\le t$,
$$
E(\langle X_{T_n}(t),\varphi\rangle-\langle X_{T_n}(s),\varphi\rangle)^2= E\theta\, A'_{T_n} +({\rm Var}\, \theta-E\theta)B'_{T_n},
$$
where 
$$
A'_T=\frac{T^2}{F_T^2}\int_s^t\int_s^tp_{T|u-v|}(0,0)dudv,\: B'_T=\frac{T^2}{F_T^2}\int_s^t\int_s^tp_{T(u+v)}(0,0)dudv.
$$
We have 
$$
A'_T\le 2T^{1+\gamma}\int_0^{t-s}\int_0^vp_{Tr}(0,0)drdv\le 2T^\gamma(t-s)\int_0^{T(t-s)}p_r(0,0)dr\le C(t-s)^{1-\gamma}
$$
by \eqref{eq:3.5a}. As $t\mapsto p_t(0,0)$ is decreasing, then
$$
B'_T\le T^{1+\gamma}\int_0^{t-s}\int_0^{t-s}p_{Tr}(0,0)drdv\le  C(t-s)^{1-\gamma}.
$$
Hence tightness follows. This completes the proof of the theorem.
\qed
\vglue.5cm
\noindent
{\bf Proof of Theorem \ref{t:2.4}}
\vglue.4cm
\noindent
Without loss of generality we may assume that the Poisson parameter of the law of $\theta$ is 1. We follow the same scheme as in the previous theorem. Again we restrict to the case $\varphi=\UNO_{\{0\}}$.

First we show convergence of covariances. A standard argument for branching Poisson system yields
$$
{\rm Cov}(\langle N_v,\varphi\rangle,\langle N_u,\varphi\rangle)=\sum_{x\in\Omega_M}\left(\varphi(x){\cal T}_{|v-u|}\varphi(x)+V\int_0^{u\wedge v}\varphi(x){\cal T}_{u+v-2r}\varphi(x)dr\right),
$$
where, for our $\varphi$, by \eqref{eq:2.1} and \eqref{eq:2.12}, after an obvious substitution we have 
\begin{equation}
\label{eq:3.19}
{\rm Cov}(\langle X_T(s),\varphi\rangle\langle X_T(t),\varphi\rangle)=I(T)+
I\!\!I(T),
\end{equation}
where 
\begin{equation}
\label{eq:3.20}
I(T)=\frac{T^2}{T^{2-\gamma}}\int_0^t\int_0^sp_{T|v-u|}(0,0)dudv,
\end{equation}

\begin{equation}
\label{eq:3.21}
I\!\!I(T)=V\frac{T^3}{T^{2-\gamma}}\int_0^t\int_0^s\int_0^{u\wedge v}p_{T(u+v-2r)}(0,0)drdudv.
\end{equation}
Assuming that $s\le t$, 
\begin{equation}
\label{eq:3.22}
I(T)\le 2T^\gamma\int_0^t\int_0^tp_{Tu}(0,0)dudv=2tT^{\gamma-1}\int_0^{Tt}p_u(0,0)du\to 0
\end{equation}
as $T\to\infty$, by  \eqref{eq:3.5c}, since $\gamma<1$.

From \eqref{eq:3.2}, for $T_n=a^{-n}$ we have 
\begin{eqnarray}
I\!\!I(T_n)&=&V\frac{M-1}{Mb^{1+\gamma}}\int_0^t\int_0^s\int_0^{u\wedge v}\sum_{j=-n}^\infty(ba^j)^{\gamma+1}e^{-ba^j(u+v-2r)}drdudv\nonumber\\
&\to&V\frac{M-1}{Mb^{1+\gamma}}\int_0^t\int_0^s\int_0^{u\wedge v}\sum_{j\in\zetita}^\infty(ba^j)^{\gamma+1}e^{-ba^j(u+v-2r)}drdudv\nonumber\\
&=&V\frac{M-1}{Mb^{1+\gamma}}E\zeta^H_s\zeta^H_t, \label{eq:3.23}
\end{eqnarray}
by \eqref{eq:2.7}, with $H={(2-\gamma)}/{2}$. Then 
\eqref{eq:3.19}-\eqref{eq:3.23} imply the desired convergence of the covariance function of $\langle X_T,\varphi\rangle$.

Thanks to (\ref{eq:3.7}) we can prove convergence of finite dimensional distributions using the central limit theorem (Lyapunov criterion). Analogously as in the proof of Theorem 2.2 this reduces to showing that
\begin{equation}
\label{eq:3.24}
S_T:=\frac{1}{F^3_T}\sum_{x\in\Omega_M}E\left(\int^T_0\langle N^x_n,\varphi\rangle du\right)^3\to 0\,\,\,{\rm as}\,\, T\to\infty,
\end{equation}
where $N^x$ is the empirical process of a branching system starting from a single particle at $x$. Define
$$v_q(x,t)=1-E{\rm exp}\left(
-q\int^t_0\langle N^x_s,\varphi\rangle ds\right).$$
A standard argument (conditioning on the first branching and applying the Feynman-Kac formula, see e.g. \cite{GLM} or \cite{BGT1}) shows that
$$v_q(x,t)=\int^t_0{\cal T}_{t-s}\left(q\varphi(1-v_q(\cdot, s))-\frac{V}{2}v^2_q(\cdot, s)\right)(x)ds.$$
Hence, by a similar reasoning as in (3.45)-(3.47) of \cite{BGT2} (see also the proof of Theorem 2.4 in \cite{BGT9}, complete version) we obtain
\begin{eqnarray}
S_T&=&\frac{1}{F^3_T}\sum_{x\in\Omega_M}\frac{\partial^3}{\partial q^3}v_q(x,T)|_{q=0}\nonumber\\
\label{eq:3.25}
&=&I_1(T)+I_2(T)+I_3(T)+I_4(T),
\end{eqnarray}
where
\begin{eqnarray*}
I_1(T)&=&\frac{6}{F^3_T}\sum_{x\in\Omega_M}\int^T_0{\cal T}_{T-u}
\left(\varphi\int^u_0{\cal T}_{u-u_1}\left(\varphi\int^{u_1}_0{\cal T}_{u_1-u_2}\varphi du_2\right)du_1\right)(x)du,
\\
I_2(T)&=&\frac{3V}{F^3_T}\sum_{x\in\Omega_M}
\int^T_0{\cal T}_{T-u}
\left(\varphi\int^u_0{\cal T}_{u-u_1}
\left(\int^{u_1}_0{\cal T}_{u_1-u_2}\varphi du_2\right)^2du_1\right)(x)du,
\\
I_3(T)&=&\frac{6V}{F^3_T}\sum_{x\in\Omega_M}\int^T_0{\cal T}_{T-u}
\left(\left(
\int^u_0{\cal T}_{u-u_1}\varphi du_1\right)\int^u_0{\cal T}_{u-u_2}
\left(\varphi\int^{u_2}_0{\cal T}_{u_2-u_3}\varphi du_3\right)du_2\right)(x)du,
\\
I_4(T)&=&\frac{3V^2}{F^3_T}\sum_{x\in\Omega_M}\int^T_0{\cal T}_{T-u}
\left(\left(
\int^u_0{\cal T}_{u-u_1}\varphi du_1\right)\int^u_0{\cal T}_{u-u_2}
\left(\int^{u_2}_0{\cal T}_{u_2-u_3}\varphi du_3\right)^2du_2\right)(x)du.
\end{eqnarray*}
Note that in all the integrals ${\cal T}_{T-u}$ can be omitted since the counting measure is invariant for ${\cal T}_t$.

For $\varphi=\UNO_{\{0\}}$, using \eqref{eq:3.5} and then the estimate \eqref{eq:3.5c} we have
\begin{eqnarray*}
I_1(T)&\leq&\frac{C}{T^{3(2-\gamma)/2}}\int^T_0G_T(\varphi G_T\varphi)(0)du\\
&\leq&\frac{C_1}{T^{2-3\gamma/2}}\to 0.
\end{eqnarray*}
Similarly,
$$I_2(T)\leq\frac{C}{F^3_T}T\sum_{y\in\Omega_M}
\int^T_0p_r(0,y)(G_T\varphi(y))^2dr<\frac{C_1}{T^{2-3\gamma/2}}\sum_{y\in\Omega_M}(G_T\varphi(y))^2.$$
Notice that
\begin{eqnarray}
\sum_{y\in\Omega_M}(G_T\varphi(y))^2&=&\int^T_0\int^T_0\sum_{y\in\Omega_M}p_s(0,y)p_t(0,y)dsdt\nonumber\\
\label{eq:3.26}
&=&\int^T_0\int^T_0p_{s+t}(0,0)dsdt\leq CT^{1-\gamma}
\end{eqnarray}
by (\ref{eq:3.4}). Hence $I_2(T)\to 0$.

Also $I_3(T)\to 0$ by \eqref{eq:3.5c} and (\ref{eq:3.26}).

Finally, applying consecutively the Schwarz and Young inequalites (for the convolution on $\Omega_M$) and (\ref{eq:3.26}) we have
\begin{eqnarray}
I_4(T)&\leq&
\frac{T}{F^3_T}\sum_{x\in\Omega_M}G_T\varphi(x)G_T(G_T\varphi)^2(x)
\label{eq:3.26a}
\\
&\leq&\frac{T}{F^3_T}
\sqrt{\sum_{x\in\Omega_M}(G_T\varphi(x))^2}
\sqrt{\sum_{x\in\Omega_M}G_T(G_T\varphi)^2(x)}\nonumber\\
&\leq&C\frac{T^{(3-\gamma)/2}}{F^3_T}
\sqrt{\sum_{x\in\Omega_M}(G_T
\UNO_{\{0\}}(x))^2}\sum_{x\in\Omega_M}(G_T\varphi(x))^2\nonumber\\
&\leq&C_1 T^{-\gamma/2}\to 0.
\nonumber
\end{eqnarray}
Hence (\ref{eq:3.24}) follows by (\ref{eq:3.25}).

It remains to prove tightness. From (\ref{eq:3.19})-(\ref{eq:3.21}), for
$s<t$ we obtain
\begin{equation}
\label{eq:3.27}
E(\langle X_T(t),\varphi\rangle -\langle X_T(s),\varphi\rangle)^2=J_1(T)+J_2(T)
\end{equation}
where
\begin{eqnarray} 
\label{eq:3.28}
J_1(T)&=&T^\gamma\int^t_s\int^t_sp_{T|v-u|}(0,0)dudv,\\
\label{eq:3.29}
J_2(T)&=&VT^{1+\gamma}\int^t_s\int^t_s\int^{u\wedge v}_0p_{T(u+v-2r)}(0,0)drdudv.
\end{eqnarray}

We have
\begin{eqnarray}
J_1(T)&\leq&2T^\gamma(t-s)\int^{t-s}_0p_{Tv}(0,0)dv\nonumber\\
&=&2T^{\gamma-1}(t-s)\left(\int^{T(t-s)}_0p_v(0,0)dv\right)^{1-\gamma}
\left(\int^{T(t-s)}_0p_v(0,0)dv\right)^\gamma\nonumber\\
\label{eq:3.30}
&\leq&C(t-s)^{2-\gamma},
\end{eqnarray}
by \eqref{eq:3.5c} and since $p_v(0,0)\leq 1$.

Next, by (\ref{eq:3.4}),
\begin{eqnarray}
\label{eq:3.31}
J_2(T)&\leq&C\int^t_s\int^t_s\int^{u\wedge v}_0|u+v-2r|^{-\gamma-1}drdudv\\
\label{eq:3.32}
&\leq&C_1\int^t_s\int^t_s|u-v|^{-\gamma}dudv=C_2(t-s)^{2-\gamma}.
\end{eqnarray}
Combining (\ref{eq:3.27})-(\ref{eq:3.32}) we obtain tightness.
\hfill$\Box$
\vglue.5cm
\noindent
{\bf Proof of Theorem 2.6} 
\vglue.5cm
\noindent
The proof is a modification of the proof of Theorem 2.4. Again, we have 
(\ref{eq:3.19})-(\ref{eq:3.21}) since $H_T$ cancels out (see (\ref{eq:2.15})). The limit of covariances is obtained in the form (\ref{eq:2.7}) by proving 
(\ref{eq:3.22}) and (\ref{eq:3.23}). The only difference is that in (\ref{eq:3.22}) we use \eqref{eq:3.5a} instead of \eqref{eq:3.5c}. To prove convergence of finite dimensional distributions we show (\ref{eq:3.24}). $S_T$ has the form (\ref{eq:3.25}) with $H_TI_i(T)$ instead of $I_i(T), i=1,\ldots,4$. Convergence to $0$ of 
$H_TI_i(T), i\leq 3$, is easily obtained from \eqref{eq:3.5a} and (\ref{eq:3.26}) without any condition on the rate of divergence of $H_T$. For $H_TI_4(T)$ we use \eqref{eq:3.26a}, hence, by \eqref{eq:3.5a},
\begin{eqnarray*}
H_TI_4(T)&\leq&C\frac{T^{1-\gamma}H_T}{F^3_T}\sum_{x\in\Omega_M}\int^T_0p_t(0,x)dt\sum_y\int^T_0p_s(x,y)ds\int^T_0p_u(0,y)du\\
&=&C\frac{1}{H^{1/2}_TT^{2-\gamma/2}}\int^T_0\int^T_0\int^T_0p_{t+s+u}(0,0)dtdsdu\\
&\leq&C_1T^{-\gamma/2}H^{-1/2}_T\to 0,
\end{eqnarray*}
by (\ref{eq:2.14}), where in the last estimate we have used (\ref{eq:3.4}). Hence (\ref{eq:3.24}) follows.

To prove tightness on $[0,\tau]$ for any fixed $\tau$ we go to 
(\ref{eq:3.27})-\eqref{eq:3.29} ($H_T$ cancels out). As $\gamma<0$, we have
$$J_1(T)\leq C(t-s)^2,$$
and from (\ref{eq:3.31}) it follows that
$$J_2(T)\leq C\int^t_s\int^t_s(u+v)^{-\gamma}dudv\leq C(\tau)(t-s)^2.$$
\hfill$\Box$
\vglue.5cm
\noindent
{\bf Proof of Proposition 2.8}
\vglue.5cm
\noindent
Modifying the argument at the beginning of the proof of Theorem 2.2 we have for $\varphi,\psi\in {\cal B}_b$ and $s\leq t$,
\begin{eqnarray*}
\lefteqn{{\rm Cov}\left(\langle X_T(s),\varphi\rangle, \langle X_T(t),\psi\rangle\right)=E\theta\frac{T^2}{F^2_T}\sum_{x\in\Omega_M}\int^s_0\int^t_0\varphi(x)
{\cal T}_{T|u-v|}\psi(x)dudv}\\
&=&E\theta\sum_{x\in\Omega_M}\left(\int^s_0\int^{T(t-v)}_0\varphi(x){\cal T}_u\psi(x)dudv+\int^s_0\int^{T(s-u)}_0\varphi(x){\cal T}_v\psi(x)dvdu\right)\\
&\to&2E\theta s\sum_{x\in\Omega_M}\varphi(x)G\psi(x)
\end{eqnarray*}
as $T\to\infty$.

To prove the convergence in law of the sum $\sum ^m_{k=1}\langle X_T(t_k),\varphi_k\rangle, \varphi_1,\ldots,\varphi_k\in {\cal B}_b,t_1,\ldots,t_k\geq 0$, we use the Lyapunov criterion, which again reduces to showing  (\ref{eq:3.18}) with $T$ instead of $T_n$ and $\varphi\geq 0$. Formula \eqref{eq:3.18a} implies that
$$S_T\leq\frac{6}{F^3_T}T\sum_{x\in\Omega_M}\varphi(x)G(\varphi G\varphi)(x)\to 0.$$
\hfill$\Box$
\vglue.5cm
\noindent
{\bf Proof of Proposition 2.10}
\vglue.5cm
\noindent
(a) follows from (\ref{eq:2.3}) and the known formulas
\begin{eqnarray}
\label{eq:3.33}
t&=&\frac{2}{\pi}\int^\infty_0\frac{{\rm sin}^2(tu)}{u^2}du,\quad t>0,\\
\label{eq:3.34}
1-e^{-rt}&=&\frac{2}{\pi}r\int^\infty_0
\frac{{\rm sin}^2(tu)}{{r^2}/{4}+u^2}du,\quad t>0.
\end{eqnarray}
(b) is a direct consequence of (a) and (\ref{eq:1.3}).

\noindent
(c) Formula (\ref{eq:2.21}) for osfBm also follows from (a) and (\ref{eq:1.4}). To derive it for onsfBm we interchange the sum and the integrals in 
(\ref{eq:2.7}), and after integrating we use (\ref{eq:3.33}) and (\ref{eq:3.34}).
\vglue.5cm
\noindent
{\bf Comments on the proofs of Propositions 2.12, 2.14, 2.16, 2.17 and 2.20}
\vglue.5cm
\noindent
In Proposition 2.12, (\ref{eq:2.22})-(\ref{eq:2.24}) follow by direct calculations (for osfBm and onsfBm \eqref{eq:2.7}  is again useful here). Formula 
(\ref{eq:2.22}) has already appeared in (\ref{eq:3.15}).

Proposition 2.14 can also be verified directly (convergence in $L^2$ and hence almost surely as well, by independence).

In Proposition 2.16 formula (\ref{eq:2.27}) is also immediate once we know that the process 
$\varrho^U$ exists. 
$\tilde{\tilde{h}}^U_t$ defined by 
(\ref{eq:2.29}) is finite for $U\in (0,1/2)$ by (\ref{eq:3.3}). Positive definiteness of (\ref{eq:2.28}) follows from
$$\frac{1}{ba^j}
(e^{-ba^j|s-t|}-
e^{-ba^j(s+t)})=
2\int^\infty_0
\UNO_{[0,s]}(r)e^{-(s-r)ba^j}
\UNO_{[0,t]}(r)e^{(t-r)ba^j}dr.$$

In Proposition 2.17 the estimates (\ref{eq:2.30}) for ofBm follow from 
(\ref{eq:2.4}). For osfBm and onsfBm we use (\ref{eq:2.7}) and (\ref{eq:3.3}), and then (\ref{eq:2.30}) 
is a consequence of the corresponding estimates for sfBm and nsfBm (see (2.4) and (2.5) in \cite{BGT0} and Proposition 3.4 \cite{BGT5}).

Finally, Proposition 2.20 can be also reduced to the corresponding properties of the non-oscillatory processes employing (\ref{eq:3.15}), (\ref{eq:2.7}) and
(\ref{eq:3.3}).

\vglue.5cm
\noindent
{\bf Proof of Proposition 2.18}
\vglue.5cm
\noindent
(a) First we consider ofBm $\vartheta=\xi^H$. By stationarity of increments we have $EV_\varepsilon(\xi^H)/h_\varepsilon=T$, hence for the convergence in $L^2$ in \eqref{eq:2.30b} it suffices to prove that
\begin{equation}
\label{eq:3.34a}
\lim_{\varepsilon\to 0}{\rm Var}
\left(\frac{1}{h_\varepsilon}V_\varepsilon(\xi^H)\right)
=0.
\end{equation}
Using \eqref{eq:2.30a} and the identity
\begin{equation}
\label{eq:3.34b}
EX^2Y^2=EX^2EY^2+2(EXY)^2,
\end{equation}
which holds for any centered Gaussian vector $(X,Y)$, we obtain
\begin{equation}
\label{eq:3.34c}
{\rm Var}\left(\frac{1}{h_\varepsilon}V_\varepsilon(\xi^H)\right)=2\int^T_0\int^T_0(I_\varepsilon(s,t))^2dsdt,
\end{equation}
where
\begin{equation}
\label{eq:3.34d}
I_\varepsilon(s,t)=\frac{1}{\varepsilon^{2H}h_\varepsilon}E(\xi^H_{s+\varepsilon}-\xi^H_s)(\xi^H_{t+\varepsilon}-\xi^H_t)\geq 0.
\end{equation}
By the Schwarz inequality,
\begin{equation}
\label{eq:3.34e}
I_\varepsilon(s,t)\leq 1.
\end{equation}
 On the other hand, for $s<t$, from \eqref{eq:2.22} and
 \eqref{eq:3.3} we have
$$I_\varepsilon(s,t)\leq\frac{1}{\varepsilon^{2H}}\int^{t+\varepsilon}_t\int^{s+\varepsilon}_s|u-v|^{2H-2}dudv\to 0$$
since $H\in(1/2,1)$. Hence \eqref{eq:3.34a} follows.

For the proof of almost sure convergence we need a finer estimate of the variance. It is easy to see that for $t-s>2\varepsilon$ we have
$$I_\varepsilon(s,t)\leq C\varepsilon^{2-2H}|t-s|^{2H-2},$$
hence, by \eqref{eq:3.34e},
$$(I_\varepsilon(s,t))^2\leq C_1
\left(\frac{\varepsilon}{|t-s|}\right)^r\,\, \hbox{\rm for any}\,\, 0\leq r\leq 4-4H.$$
Therefore, from \eqref{eq:3.34c}  we obtain
\begin{equation}
\label{eq:3.34f}
{\rm Var}\left(\frac{1}{h_\varepsilon}V_\varepsilon(\xi^H)\right)\leq C_2T\varepsilon+C_3(T)\varepsilon^r\leq C(T)\varepsilon^r
\end{equation}
for any $0\leq r<(4-4H)\wedge 1$.

Now we take $\varepsilon_n$ and $\delta$ as in the Proposition and fix $\delta <r< (4-4H)\wedge 1$. Using
$$P\left(\left|\frac{1}{h_{\varepsilon_n}}V_{\varepsilon_n}(\xi^H)-T\right|>
\varepsilon_n^{{(r-\delta)}/{2}}\right)\leq \varepsilon_n^{-(r-\delta)}
{\rm Var}\left(\frac{1}{h_{\varepsilon_n}}V_{\varepsilon_n}(\xi^H)\right),$$
\eqref{eq:3.34f} and the Borel-Cantelli lemma, the assertion follows.

Next we prove \eqref{eq:2.30b} for osfBm $\vartheta=\zeta^H$. Let $\eta^H$ be a continuous centered Gaussian process with covariance \eqref{eq:2.8} independent of $\zeta^H$. We know that $\zeta^H+\eta^H/\sqrt{2}$ is ofBm. Hence it is easy to see that by \eqref{eq:2.30b} for ofBm, the proof reduces to showing that
$$\lim_{\varepsilon\to 0}V_\varepsilon(\eta^H)=0$$
in $L^2$ and almost surely. The $L^2$ convergence follows easily from \eqref{eq:3.34b} and the estimate
$$|E(\eta^H_{s+\varepsilon}-\eta^H_s)(\eta^H_{t+\varepsilon}-\eta^H_t)|\leq C
\int^{t+\varepsilon}_t\int^{t+\varepsilon}_s(u+v)^{2H-2}dudv,$$
which is a consequence of \eqref{eq:3.16} and \eqref{eq:3.3}. Recall that the expression in \eqref{eq:3.16} is, up to a constant,
$$E(\eta^H_s\eta^H_t),\,\, H=(1-\gamma)/2.$$

In the proof of almost sure convergence we use the H\"older continuity of the paths of $\eta^H$ with any exponent $q<H$, and the representation
$$\eta^H_t=\int^t_0r_udu$$
for a Gaussian process $r_u$ (see Subsection 2.3.5). We then have
\begin{eqnarray*}
V_\varepsilon(\eta^H)&\leq& C\varepsilon^{q-2H}\int^T_0\int^{t+\varepsilon}_t|r_u|dudt\\
&\leq&C_1\varepsilon^{q-2H-1}\to 0
\end{eqnarray*}
for $q$ sufficiently close to $H$.

\noindent
(b) The proof is similar to  that of part (a), and we skip it.

\hfill$\Box$
\vglue.5cm
\noindent
{\bf Proof of Theorem 2.21}
\vglue.5cm
\noindent
We have
$$G\varphi(x)=\sum_{y\in\Omega_M}G(x-y)\varphi(y),$$
where
\begin{equation}
\label{eq:3.35}
G(y)=\left\{\begin{array}{lll}
D&{\rm if}&y=0,\\
{A}/{c^{i-1}}&{\rm if}& |y|=i\geq 1,
\end{array}\right.
\end{equation}
with $A=(M-c)/(Mb(c-1))$ 
and $D$ given by (\ref{eq:2.31}) (see  \cite{DGW1}, it can also be easily obtained from (\ref{eq:3.2})).

By (\ref{eq:2.17}), for $s\leq t$,
\begin{equation}
\label{eq:3.36}
E\widetilde{\vartheta}_s\widetilde{\vartheta}_t=2\langle G\UNO_{B_t},
\UNO_{B_s}\rangle=2(I+I\!\!I),
\end{equation}
where
\begin{eqnarray}
\label{eq:3.37}
I&=&\sum_{x\in B_s}\sum_{y\in B_s}G(x-y),\\
\label{eq:3.38}
I\!\!I&=&\sum_{x\in B_s}\sum_{y\in B_t\backslash B_s}G(x-y).
\end{eqnarray}

From the hierarchical structure of $\Omega_M$ it follows that if $x\in B_s$, then $y\in B_s$ if and only if $x-y\in B_s$. Let $S_i=\{x\in\Omega_M:|x|=i\}$. The cardinalities of $S_i$ and $B_i$ are $M^{i-1}(M-1), i=1,2,\ldots,$ and 
$M^i, i=0,1,\ldots$, respectively.
Hence
\begin{eqnarray}
I&=&\sum_{x\in B_s}\sum_{y\in B_s}G(y)=M^s\sum^s_{i=0}\sum_{y\in S_i}G(y)\nonumber\\
\label{eq:3.39}
&=&M^s\left(D+\sum^s_{i=1}(M-1)M^{i-1}\frac{A}{c^{i-1}}\right)=DM^s
\left(\frac{M}{c}\right)^s,
\end{eqnarray}
by (\ref{eq:3.35}).

By the hierarchical structure of $\Omega_M$, if $|y|=s+k$ and $|x|\leq s$, then $|x-y|=s+k$. Therefore, for $t>s$,
\begin{eqnarray}
I\!\!I&=&\sum^{t-s}_{k=1}\sum_{y\in S_{s+k}}\sum_{x\in B_s}G(x-y)\nonumber\\
&=&\sum^{t-s}_{k=1}(M-1)M^{s+k-1}M^s\frac{A}{c^{s+k-1}}\nonumber\\
\label{eq:3.40}
&=&DM^s\left(\left(\frac{M}{c}\right)^t-\left(\frac{M}{c}\right)^s\right).
\end{eqnarray}
Then (\ref{eq:3.36})-(\ref{eq:3.40}) imply (\ref{eq:2.33}).
\hfill$\Box$

\section{Comments}

In this section we discuss some possible extensions of the results 
of Section 2.

\subsection{Oscillatory processes with small parameters}
\label{sub:4.1}

We have defined ofBm and osfBm for 
$H\in ({1}/{2},1)$ only, whereas fBm and sfBm are well defined for 
$H\in (0,1)$. In the oscillatory case $h$ given by 
(\ref{eq:2.3}) is infinite for $H\in (0,{1}/{2})$, so a direct extension is impossible. It seems that a natural way to define $\xi^H$ and $\zeta^H$ for 
$H\in (0,{1}/{2})$ 
is to consider the covariances (\ref{eq:1.3}) and (\ref{eq:1.4}) with 
$\tilde{\tilde h}^H$ given by (\ref{eq:2.29}). These processes can be extended to the whole real line (see Remark 2.11(b)). We then have
$$(\zeta^H_t)_{t\geq 0}\stackrel{d}{=}
\left(\frac{\xi^H_t+\xi^H_{-t}}{\sqrt{2}}\right)_{t\geq 0},\,\,(\varrho^H_t)_{t\geq 0}
\stackrel{d}{=}
\left(\frac{\xi^H_t-\xi^H_{-t}}{\sqrt{2}}\right)_{t\geq 0},$$
where $\varrho^H$ is the process that appears in Proposition 2.16. In other words, $\int^t_0\varrho^H_udu$ is an onsfBm with parameter $H+1$, so our definitions of $\xi^H,\zeta^H$ and $\varrho^H$ yield the same relationships as in the non-oscillatory case (cf. \cite{DZ,BGT0,BGT5}). Existence of $\xi^H$ (positive definiteness of the covariance) follows from \eqref{eq:3.34}. It is not clear what would be particle picture interpretations of the oscillatory processes $\xi^H, \zeta^H$ and $\varrho^H$ with $H\in (0,{1}/{2})$, as the method used in \cite{BT} for the non-oscillatory case is not applicable here.

\subsection{Models with immigration}
\label{sub:4.2}

Similarly to the case of $\alpha$-stable motions in ${\mathbb R}^d$, one can add immigration to the particle systems as in \cite{GNR}. We assume that particles immigrate according to a space-time uniform Poisson random field on 
$\Omega_M\times{\mathbb R}_+$. The new particles evolve independently, undergoing $c$-random walks with or without branching. In \cite{GNR} limits of covariances of occupation time fluctuations were calculated, thus some new covariance functions were obtained, and the corresponding Gaussian processes were long range dependent. Analogous calculations can be carried out for the hierarchical case.

In the non-branching case, for $\gamma<0\,\, (c<1), T_n=a^{-n}$ and 
$F_n=T_n^{1-{\gamma}/{2}}$, taking for simplicity 
$\varphi=\UNO_{\{0\}}$, it can be shown that
\begin{equation}
\label{eq:4.1}
\lim_{n\to\infty}{\rm Cov}(\langle X_{T_n}(s),\varphi\rangle, 
\langle X_{T_n}(t),\varphi\rangle)
= K\int^{s\wedge t}_0 u\left((t-u)^{-\gamma}
\tilde{\tilde h}_{t-u}+(s-u)^{-\gamma}\tilde{\tilde h}_{s-u}\right)du,
\end{equation}
where $\tilde{\tilde h}$ is given by (\ref{eq:2.29}) with 
$U=-{\gamma}/{2}$. The right-hand side of (\ref{eq:4.1}) is an oscillatory 
analogue of the covariance of weighted fractional Brownian motion with parameters 1 and $-\gamma$ as defined in \cite{BGT5}. Another form of the covariance (\ref{eq:4.1}) is
$$\int^{s\wedge t}_{0}E\xi^H_{s-u}\xi^H_{t-u}du,$$
where $\xi^H$ is ofBm with $H={(1-\gamma})/{2}$.

In the branching case, for $0<\gamma<1\,\, (1<c<\sqrt{M}), T_n=a^{-n}, F_n=
T_n^{{(3-\gamma)}/{2}},\varphi=\UNO_{\{0\}}$, it can be shown that the limit of covariances is, up to a constant, the difference of the covariance (\ref{eq:4.1}) with $\gamma-1$ instead of $\gamma$, and the covariance of onsfBm with 
$H={(3-\gamma)}/{2}$.

We have not attempted to prove convergence of finite dimensional distributions of $X_{T_n}$.

\subsection{Other initial configurations}
\label{sub:4.3}

One can consider initial conditions different from those described in Subsection 2.1. For example, it is possible to study a branching system starting from an 
equilibrium measure which is known to exist for $\gamma>0$ (see Appendix in \cite{DGW1}, and \cite{GW2}). By analogy with the result for the $\alpha$-stable case in $\mathbb R^d$ obtained in \cite{M1}, one should expect that in the hierarchical model with $0<\gamma<1$ the limit process is ofBm. Preliminary calculations confirm this conjecture.

Another possibility is to consider inhomogeneous initial configurations analogous to those studied in \cite{BGT6} and \cite{BGT8}.

\subsection{More general random walks}
\label{sub:4.4}

It is possible to consider general $r_j$-random walks as described in the Introduction, i.e. $(r_j)_{j=1,2,\ldots}$ is a probability distribution 
not necessarily of the form (\ref{eq:1.7}). It seems that, similarly as in \cite{BGT10}, a natural assumption in this case is the existence of a positive limit
$$\lim_{j\to\infty}\frac{r_{j+1}}{r_j}=a.$$
We conjecture that if this limit is strictly less than $1$, then it plays the role of $a$ (see 
(\ref{eq:2.2})) and the results should be similar to those obtained for $c$-random walks. If $a=1$, then the limits should be as in Proposition 2.8 or Remark 2.9(a). The degree is given by  $\gamma=\log_{1/a}M-1$ if $a<1$, and $\gamma=\infty$ if $a=1$ (see the Appendix of \cite{BGT10}, or \cite{DGW4}). The case of ``small'' $\gamma$ seems more difficult.

\subsection{Infinite variance branching}
\label{sub:4.5}

Analogously as in \cite{BGT3,BGT4}, instead of binary branching one can consider a branching mechanism with probability generating function
$$s+\frac{(1-s)^{1+\beta}}{1+\beta},\quad 0<s<1,$$
where $0<\beta<1$. The assumption $0<\gamma<1$ in Theorem 2.4 is replaced by 
${1}/{\beta}-1<\gamma<{1}/{\beta}$, and $F_t=T^{{(2-\beta\gamma)}/{(1+\beta)}}$. Some preliminary calculations yield $(1+\beta)$-stable limits of rather complicated forms. Again they have dependent increments.

Note that the condition ${1}/{\beta}-1<\gamma$ implies persistence of the system, since it can be shown that
$$\sum_{x\in\Omega_M}\int^\infty_0({\cal T}_t\varphi(x))^{1+\beta}dt<\infty,\quad \varphi\in{\cal B}_b, \quad \varphi\geq 0,$$
 which is sufficient for persistence (see the Appendix of \cite{DGW1}). Persistence means that the empirical process $N_t$ has a limit in distribution as $t\to\infty, N_\infty$, with full intensity measure, i.e., $EN_\infty=EN_0$, assuming that $N_0$ is homogeneous Poisson (see \cite{GW1,GW2} and references therein).

\subsection{Oscillatory bi-fractional Brownian notion}
\label{sub:4.6}

Bi-fractional Brownian motion  has covariance function
$$\frac{1}{2^K}\Big((s^{2H}+t^{2H})^K-|t-s|^{2HK}\Big),$$
with parameters $H\in (0,1], K\in (0,1]$. It is an extension of fBm which was obtained analytically \cite{HV}. In the same way as in \cite{HV}, one can obtain an oscillatory analogue of ofBm with covariance function
$$\frac{1}{2^K}\Big((s^{2H}h_{s}+t^{2H}h_{t})^K-(|s-t|^{2H}h_{|s-t|})^K\Big),$$%
$h$ defined by \eqref{eq:2.3}, $H\in (1/2,1], K\in (0,1)$.
We have not found particle picture interpretations for bi-fractional Brownian motion and its oscillatory analogue.

\noindent
\def\refname{\hbox{\normalsize\bf References}}

\end{document}